\documentclass[sigconf]{acmart}

\AtBeginDocument{%
  }

\copyrightyear{2024}
\acmYear{2024}
\setcopyright{rightsretained}
\acmConference[KDD '24]{Proceedings of the 30th ACM SIGKDD Conference on Knowledge Discovery and Data Mining}{August 25--29, 2024}{Barcelona, Spain}
\acmBooktitle{Proceedings of the 30th ACM SIGKDD Conference on Knowledge Discovery and Data Mining (KDD '24), August 25--29, 2024, Barcelona, Spain}\acmDOI{10.1145/3637528.3671783}
\acmISBN{979-8-4007-0490-1/24/08}

\acmSubmissionID{rtp0626}

\usepackage{algorithm}
\usepackage{algorithmic}
\usepackage{hyperref}
\usepackage{url}
\usepackage{minitoc}
\usepackage{booktabs}
\usepackage{graphicx} 
\usepackage{float} 
\usepackage{subfigure} 
\usepackage{multirow}
\usepackage{enumitem}
\usepackage{bm}

\newcommand{\bzi}{\mathbf{z}_{i}}
\newcommand{\bzx}{\mathbf{z}_{\mathbf{x}}}
\newcommand{\bzxt}{\mathbf{z}_{\mathbf{x}}^{(t)}}
\newcommand{\bzxzero}{\mathbf{z}_{\mathbf{x}}^{(0)}}
\newcommand{\bzxone}{\mathbf{z}_{\mathbf{x}}^{(1)}}
\newcommand{\bzxtwo}{\mathbf{z}_{\mathbf{x}}^{(2)}}
\newcommand{\bzxT}{\mathbf{z}_{\mathbf{x}}^{(T)}}
\newcommand{\bzxpt}{\mathbf{z}_{\mathbf{x}}^{(t-1)}}
\newcommand{\bI}{\mathbf{I}}
\newcommand{\bepsilon}{\bm \epsilon}

\newcommand{\bmu}{\bm \mu}
\newcommand{\bSigma}{\mathbf{\Sigma}}
\newcommand{\ff}{\mathbf{f}}
\newcommand{\dd}{\mathbf{d}}
\newcommand{\xx}{\mathbf{x}}
\newcommand{\hatbzx}{\hat{\mathbf{z}}_{\mathbf{x}}}

\makeatletter
\gdef\@copyrightpermission{
   \begin{minipage}{0.3\columnwidth}
     \href{https://creativecommons.org/licenses/by-nd/4.0/}{\includegraphics[width=0.90\textwidth]{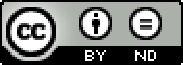}}
   \end{minipage}\hfill
   \begin{minipage}{0.7\columnwidth}
     \href{https://creativecommons.org/licenses/by-nd/4.0/}{This work is licensed under a Creative Commons Attribution-NoDerivs International 4.0 License.}
   \end{minipage}
   \vspace{5pt}
}
\makeatother
\begin{document}

\title{Effective Generation of Feasible Solutions for Integer Programming via Guided Diffusion }


\author{Hao Zeng}
\email{zenghao.zeng@cainiao.com}
\affiliation{%
  \institution{Cainiao Network}
  \city{Hangzhou}
  \country{China}
}

\author{Jiaqi Wang}
\email{tangqiao.wjq@cainiao.com}
\affiliation{%
  \institution{Cainiao Network}
  \city{Hangzhou}
  \country{China}
}

\author{Avirup Das}
\email{avirup.das@postgrad.manchester.ac.uk}
\affiliation{%
 \institution{University of Manchester}
 \city{Manchester}
 \country{United Kingdom}
 }

\author{Junying He}
\email{junying.hjy@cainiao.com}
\affiliation{%
  \institution{Cainiao Network}
  \city{Hangzhou}
  \country{China}
  }

\author{Kunpeng Han}
\email{kunpeng.hkp@cainiao.com}
%
\author{Haoyuan Hu}
\email{haoyuan.huhy@cainiao.com}
\affiliation{%
  \institution{Cainiao Network}
  \city{Hangzhou}
  \country{China}
}

\author{Mingfei Sun}
\authornote{Corresponding author.}
\email{mingfei.sun@manchester.ac.uk}
\affiliation{%
  \institution{University of Manchester}
  \city{Manchester}
  \country{United Kingdom}
}

\renewcommand{\shortauthors}{Hao Zeng et al.}

\begin{abstract}
Feasible solutions are crucial for Integer Programming (IP) since they can substantially speed up the solving process. In many applications, similar IP instances often exhibit similar structures and shared solution distributions, which can be potentially modeled by deep learning methods. Unfortunately, existing deep-learning-based algorithms, such as Neural Diving \citep{nair2020solving} and Predict-and-search framework \citep{han2023a}, are limited to generating only partial feasible solutions, and they must rely on solvers like SCIP and Gurobi to complete the solutions for a given IP problem. In this paper, we propose a novel framework that generates \emph{complete} feasible solutions \emph{end-to-end}. Our framework leverages contrastive learning to characterize the relationship between IP instances and solutions, and learns latent embeddings for both IP instances and their solutions. Further, the framework employs diffusion models to learn the distribution of solution embeddings conditioned on IP representations, with a dedicated guided sampling strategy that accounts for both constraints and objectives. We empirically evaluate our framework on four typical datasets of IP problems, and show that it effectively generates complete feasible solutions with a high probability (> 89.7 \%) without the reliance of Solvers and the quality of solutions is comparable to the best heuristic solutions from Gurobi. Furthermore, by integrating our method's sampled partial solutions with the CompleteSol heuristic from SCIP \citep{maher2017scip}, the resulting feasible solutions outperform those from state-of-the-art methods across all datasets,  exhibiting a 3.7 to 33.7\% improvement in the gap to optimal values, and maintaining a feasible ratio of over 99.7\% for all datasets.
\end{abstract}


\maketitle

\section{Introduction}
Integer Programming (IP) in the field of operation research is a class of optimization problems 
where some or all of the decision variables are constrained to be integers~\citep{wolsey1998}. 
Despite their importance in a wide range of applications such as production planning~\citep{silver1998inventory, pochet2006}, resource allocation~\citep{katoh1998resource}, and scheduling~\citep{toth2002vehicle,pantelides1995short,sawik2011}, 
IP is known to be NP-hard and in general very difficult to solve. 
For decades, a significant effort has been made to develop sophisticated algorithms and efficient solvers, 
e.g., branch-and-bound~\citep{lawler1966branch}, cutting plane method~\citep{kelley1960cutting} and large neighborhood search algorithms~\citep{pisinger2019large}. 
These methods, however, can be computationally expensive because the search space for large-scale problems can be exponentially large. 
Moreover, these algorithms rely heavily on a feasible solution input that will crucially determine the whole search process. Consequently, existing solvers, such as SCIP \citep{maher2017scip} and Gurobi \citep{gurobi2021gurobi}, firstly employ heuristic algorithms to identify high-quality feasible solutions that serve as initial starting points for the optimization process. However, the heuristic algorithms usually fail to capture the similar structure among different IP instances and the quality of initial solutions is usually low. 
Hence, having a data-driven method that produces high-quality feasible solutions for any IP instances is desirable for many real-world applications. 

To generate feasible solutions, prior works \citep{nair2020solving, han2023a, yoon2022confidence} have employed the advantage of deep learning to capture similarity of the IP instances from the same domain in order to expedite solving. For presenting IP instances via neural network, Neural Diving \citep{nair2020solving, yoon2022confidence} adopt the methodology delineated by \citet{gasse2019exact}, which models the formulation of IP instances as bipartite graphs and subsequently leveraging Graph Neural Networks (GNN) to derive variable features from these graph representations. Subsequently, a solution prediction task is employed to learn the relationship between IP instances and their solutions, with the aim of directly predicting those solutions. However, it is difficult to produce complete feasible solutions as it fails to explicitly integrate objective and constraint information during the sampling process. Neural Diving thus focus on generating partial solutions by GNN, 
where only a subset of variables is assigned values using neural networks.
Importantly, in many cases, the proportion of variables predicted by the neural network is set at a relatively low ratio (less than 50\%) to ensure feasibility.
Furthermore, such methods tend to be inefficient, primarily due to the introduction of auxiliary problems for filling in the remaining variables. 
For instance, the Completesol heuristic~\citep{maher2017scip}, a classical approach, solves an auxiliary integer programming model which is constructed by adding constraints to fix the variables from partial solutions. Nonetheless, infeasibility can arise in auxiliary problems due to the potential for partial assignments to clash with the initial constraints. To deal with infeasible assignments, in another approach, \citet{han2023a} proposes a predict-and-search algorithm which constructs a trust region based on predicted partial solutions and then search for high-quality feasible solutions via a solver. In summary, these methods require the construction of auxiliary problems to obtain feasible solutions and fail to utilize the complete information from the IP instance, as only partial variables are assigned. This situation highlights the necessity for creating an end-to-end deep learning framework capable of generating complete and feasible solutions for IP problems.

\begin{figure*}[ht] 
\centering 
\includegraphics[width=0.8\linewidth]{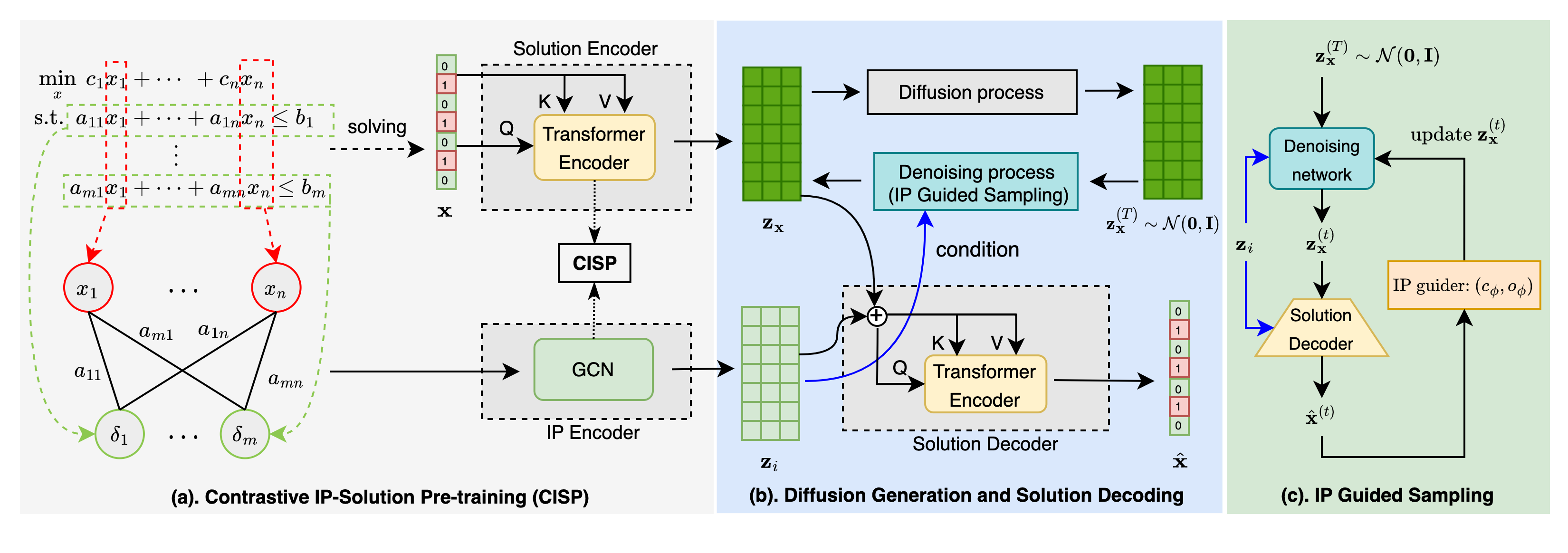}
\caption{Our method first trains the IP Encoder and Solution Encoder to acquire the IP embedding ($\bzi$) and Solution embedding ($\bzx$) using CISP. 
We then jointly train diffusion models and the solution decoder to capture the distribution of solutions given a specific IP instance. 
In the sampling stage, we employ an IP guided diffusion sampling to account for both the objective and constraints.}
\label{fig:framework}
\end{figure*}

Recently, diffusion models~\citep{ho2020denoising, sohl2015deep} have exhibited notable advantages in various generative tasks, 
primarily owing to their superior mode-coverage and diversity \citep{bayat2023a}. 
Notable applications include high-fidelity image generation \citep{dhariwal2021diffusion}, image-segmentation \citep{amit2021segdiff}, and text-to-image synthesis \citep{ramesh2022hierarchical}. 
These successes motivates the launch of an investigation into harnessing the generative capability of diffusion models for acquiring feasible solutions of IP problems.  
 
To this end, we introduce a comprehensive end-to-end generative framework presented in Figure \ref{fig:framework} to produce high-quality feasible solutions for IP problems. 
First of all, stemming inspiration from DALL.E-2~ \citep{ramesh2022hierarchical} for text-to-image translation, 
we employ a multimodal contrastive learning approach, akin to the CLIP Algorithm~\citep{radford2021learning}, 
to obtain embeddings for an IP instance $i$, denoted as $\bzi$, and solution embeddings $\bzx$ for solutions $\mathbf{x}$ (Section \ref{sec:CISP}). 
Subsequently, we employ DDPM \citep{ho2020denoising} to model the distribution of $\bzx$ conditioned on $\bzi$ (Section \ref{sec:diffusion}). 
During this phrase, a decoder is concurrently trained with the task of solution reconstruction (Section \ref{sec:decoder}). 
Finally, to enhance the quality of the feasible solutions during the sampling process, we propose the IP-guided sampling approaches tailored for both DDPM and DDIM~\citep{song2020denoising} which explicitly consider both constraints and objectives during sampling. 
Our experimental results shown in Section \ref{sec:experiments}  substantiate the efficacy of this approach in generating complete and feasible solutions for a given IP instance with a higher probability. Besides, by combining with the CompleteSol heuristic, the solutions from our methods have better quality than the state-of-the-art. Importantly, to the best of our knowledge, our approach is the first to have the ability to generate complete and feasible solutions using \emph{pure} neural techniques, without relying on any solvers.

\section{Background}
\paragraph{\textbf{Integer Programming and Its Representations.}}
Integer programming (IP) is a class of NP-hard problems where the goal is to optimize a linear objective function, 
subject to linear and integer constraints. 
Without loss of generality, we focus on minimization which can be formulated as follows,
\begin{equation}
    \label{eq:ip}
    \min_{\mathbf{x}} \mathbf{c}^\top \mathbf{x} \qquad \text{subject to } \mathbf{Ax}\leq \mathbf{b}, \qquad \mathbf{x}\in \mathbb{Z}^{n}
\end{equation}
where $\textbf{c}\in \mathbb{R}^n$ denotes the objective coefficient, $\mathbf{A} = [\mathbf{a}_1^\top, \mathbf{a}_2^\top, ..., \mathbf{a}_m^\top] \in \mathbb{R}^{m\times n}$ is the coefficient matrix of constraints 
and $\mathbf{b} = [b_1, b_2, ..., b_m]^\top\in \mathbb{R}^{m}$ represents the right-hand-side vector. 
For simplicity, we focus on binary integer variables, where $\mathbf{x}$ takes values in $\{0,1\}^n$. 
Throughout this paper, we adopt the term \emph{IP instance} to denote a specific instance within the domain of some Integer Programming (IP) problem.

Bipartite graph representation, proposed by \citet{gasse2019exact}, is a commonly used and useful way to extract features of an IP instance for machine learning purposes. 
This representation, see the left part of Figure \ref{fig:framework} (a) for an example, divides the constraints and variables into two different sets of nodes, 
and uses a Graph Convolution Network (GCN) to learn the representation of nodes. 
Recently, \citet{nair2020solving} proposed several changes to the architecture of GCN for performance improvements. 
Therefore, in this work, we use the bipartite graph structure combined with GCN to extract the embeddings of IP instances 
(see \citep{gasse2019exact, nair2020solving} for more details).

\paragraph{\textbf{DDPM and DDIM.}}
Diffusion models learn a data distribution by reversing a gradual noising process. 
In the DDPM method~\citep{ho2020denoising}, when presented with a data point sampled from an actual data distribution, denoted as $\bzxzero \sim q(\bzx)$, a diffusion model, as described in \citep{sohl2015deep, ho2020denoising}, typically involves two distinct phases. In the forward process, a sequence of Gaussian noise is incrementally added to the initial sample over a span of $T$ steps, guided by a variance schedule denoted as $\beta_1, \beta_2, \ldots, \beta_T$. This process yields a sequence of noisy samples $\bzxone, \bzxtwo, ..., \bzxT$. Subsequently, the transition for the forward process can be described as:
$q(\bzxt|\bzxpt)= \mathcal{N}(\bzxt; \sqrt{1-\beta_t}\bzxpt, \beta_t \bI)$. 
In fact, $\bzxt$ can be sampled at any time step $t$ in a closed form employing the notations $\alpha_t:=1-\beta_t$ and $\bar{\alpha} := \prod_{s=1}^t \alpha_s$, 
$\bzxt = \sqrt{\Bar{\alpha}_t}\bzxzero + \sqrt{1-\Bar{\alpha}_t}\bepsilon$, 
where $\bepsilon\sim \mathcal{N}(0,\bI)$. 
In the reverse process (denoising process), we need to model the distribution of $\bzxpt$ given $\bzxt$ as a Gaussian distribution, which implies that 
$p_{\theta}(\bzxpt|\bzxt) = \mathcal{N}\left(\bzxpt; \bmu_{\theta}(\bzxt,t), \bSigma_{\theta}(\bzxt,t)\right)$, 
where the variance $\bSigma_{\theta}(\bzxt,t)$ can be fixed to a known constant~\citep{ho2020denoising} or learned with a separate neural network~\citep{nichol2021improved}, 
while the mean can be approximately computed by adding $\bzxzero$ as a condition,
$\bmu_{\theta}(\bzxt, t) = \frac{\sqrt{\alpha_t}(1-\Bar{\alpha}_{t-1})}{1-\Bar{\alpha}_t}\bzxt + \frac{\sqrt{\Bar{\alpha}_{t-1}}\beta_t}{1-\Bar{\alpha}_t}\bzxzero$. 

To improve the efficiency of sampling of DDPM, DDIM~\citep{song2020denoising} formulates an alternative non-Markovian noising process with the same forward marginals as DDPM, 
but rewrites the probability $p_{\theta}(\bzxpt|\bzxt)$ in reverse process as a desired standard deviation $\sigma_t$.
DDIM them derives the following distribution in the reverse process,
\begin{multline}\label{eq:ddim-reverse}
    q_{\sigma}(\bzxpt|\bzxt,\bzxzero) \\
    = \mathcal{N}\left(\bzxpt; \sqrt{\Bar{\alpha}_t} \bzxzero + \sqrt{1-\Bar{\alpha}_{t-1}-\sigma_t^2} \bepsilon^{(t)}, \sigma_t^2\bI\right),
\end{multline}
where $\bepsilon^{(t)} = (\bzxt-\sqrt{\Bar{\alpha}}\bzxzero)/(1-\Bar{\alpha})$ shows the direction to $\bzxt$.

\section{Model Architecture}
Our training dataset consists of pairs $(i, \xx)$ of IP instance and their corresponding one feasible solution $\xx$. 
Given an instance $i$, let $\bzi\in \mathbb{R}^{n\times d}$ and $\bzx \in \mathbb{R}^{n\times d}$ be the embeddings of the IP instance $i$ and the solution $\mathbf{x}$ respectively, 
where $n$ is the number of variables and $d$ is the embedding dimension. 
It is worth noting that a IP instance can have multiple different feasible solutions, 
meaning that model need to learn the distribution of feasible solutions by conditioning on a given IP instance. Our methods do not directly apply a diffusion model to learn the distribution of solutions. Instead, we use an encoder to transform the solutions $\xx\in \{0,1\}^n$ from a discrete space to a continuous embedding space, e.g. $\bzx\in \mathbb{R}^{n\times d}$. We then construct a diffusion model to learn the distribution of the solution embeddings given an IP embedding $\bzi$. Finally, a decoder is trained to recover the predicted solution $\hat{\xx}$ from the embedding $\bzx$.
To effectively build the connection between the IP instance $i$ and solution $\xx$, 
we first apply Contrastive IP-Solution Pre-training (CISP) module, 
motivated by CLIP ~\citep{radford2021learning} which is used for text-to-image generation, 
to produce IP embedding $\bzi$ and solution embedding $\bzx$. 
Overall, our model consists of three key components: 
\begin{itemize}
    \item a \emph{Contrastive IP-Solution Pre-training (CISP) module} that produces IP embeddings $\bzi$ and solution embeddings $\bzx$;
    \item a \emph{diffusion module} $p(\bzx| \bzi)$ that generates solution embedding $\bzx$ conditioned on IP embedding $\bzi$; 
    \item and a \emph{decoder module} $p(\xx|\bzx, \bzi)$ that recovers solution $\xx$ from embedding $\bzx$ conditioned on IP embedding $\bzi$.
\end{itemize}
We provide more details on each module in the following sections. 

\paragraph{\textbf{Contrastive IP-Solution Pre-training.}}
\label{sec:CISP}
Previous works ~\citep{nair2020solving, han2023a} show the crucial importance of establishing the connection between the IP instances and the solutions, 
and propose to implicitly learn this connection through the task of predicting feasible solutions. 
These approaches may not exhibit strong generalization capabilities on new instances because they only utilize the collected solutions in dataset without considering feasibility explicitly during training. 
To more effectively capture this relationship, 
we propose to employ a contrastive learning task to learn representations for IP instances and embeddings for solutions by constructing feasible and infeasible solutions. 
The intuition behind is to ensure that the IP embeddings stay close to the embeddings of their feasible solutions, and away from the embeddings of the infeasible ones. 
To avoid explicitly constructing infeasible solutions, 
we proposed Contrastive IP-Solution Pre-training (CISP) algorithm
to train IP encoder and solution encoder. 
Specifically, for IP encoder, 
we extract the representation of IP instances via a bipartite graph as ~\citep{gasse2019exact}, 
and use the structure of GCNs from Neural Diving~\citep{nair2020solving} to generate all variables' embeddings as IP embeddings $\bzi$. 
For solution encoder, 
we use the encoder of the transformer to obtain the representations of each variable as solution embeddings $\bzx$. 
Both $\bzx$ and $\bzi$ have the same dimension to compute pairwise cosine similarities later. 
Since the number of variables in different IP instances may vary, 
we perform zero-padding on $\bzi$ and dummy-padding on $\xx$ (i.e. padding 2 for 0-1 integer programming) to align the dimensions.
The zero-padding for $\bzi$ is done to ensure that the cosine similarity remains unaffected.
CISP algorithm then learns to maximize the similarity between embeddings of IP and corresponding solution pairs, 
and to minimize the similarity between the embeddings of incorrect pairs, 
which is achieved through optimizing a symmetric cross-entropy loss, as detailed in Appendix~\ref{app:cisp}. 

\begin{figure}[ht]
    \centering
    \includegraphics{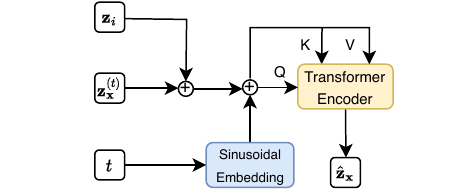}
    \caption{Diffusion model $\ff_{\theta}(\bzxt, \bzi, t)$}
    \label{fig:diffusion_model_structure}
\end{figure}

\paragraph{\textbf{Diffusion Generation.}}
\label{sec:diffusion}
To leverage diffusion models for generating feasible solutions (discrete variables), 
we use the solution embedding $\bzx\in \mathbb{R}^{n\times d}$ from the aforementioned CISP as the objective of generation. 
In addition, $\bzi$ is considered as a condition for generating high-quality results. 
According to \cite{ho2020denoising}, we parameterize 
\begin{multline}
p_{\theta}(\bzxpt|\bzxt, \bzi)\\ = \mathcal{N}\left(\bzxpt; \bmu_{\theta}(\bzxt, \bzi, t), \bSigma_{\theta}(\bzxpt, \bzi, t)\right), 
\end{multline}
$\forall t\in [T, T-1,..., 1]$ in reverse process, 
where $\bzxzero= \bzx$. 
Different from predicting the noise of each step in a general diffusion training phase, 
we predict $\bzx$ directly as it empirically performs better. 
The training loss is defined as follows,
\begin{equation}
\label{loss:diffusion}
    \mathcal{L}_{\text{MSE}} \triangleq \mathbb{E}_{t, \bzxt}\left[\|\ff_{\theta}(\bzxt,\bzi, t)-\bzx\|^2\right],
\end{equation}
where $\ff_{\theta}$ is an encoder-based transformer model with specific structure shown in Figure \ref{fig:diffusion_model_structure}  and $\bzxt  = \sqrt{\Bar{\alpha}_t}\bzx +\sqrt{1-\Bar{\alpha}_t}\bepsilon^{(t)}, \bepsilon^{(t)}\sim \mathcal{N}(\mathbf{0},\bI)$. 

\paragraph{\textbf{Solution Decoding.}}
\label{sec:decoder}
The decoder $\dd_{\phi}$ plays a crucial role in reconstructing the solution $\xx$ from the solution embedding $\bzx$. 
To enhance the robustness of the solution recovery, 
we jointly train the decoder $\dd_{\phi}$ with the diffusion model. 
Specifically, we concatenate the solution embedding $\hatbzx = \ff_{\theta}(\bzxt, \bzi, t)$ generated by the diffusion model
with the IP embedding $\bzi$, and use the concatenated vector as input to a transformer encoder to obtain the reconstructed solution $\hat{\xx} = \dd_{\phi}(\hatbzx, \bzi)$. 
This process is associated with the cross-entropy loss defined as:
$\mathcal{L}_{\text{CE}} \triangleq -\mathbb{E}_{\xx}[\log \hat{\xx}]= -\mathbb{E}_{\xx} \left[\log \dd_{\phi}(\hatbzx, \bzi) \right]$.
To explicitly account for constraints in the training process, 
we introduce a penalty term to measure the degree of constraint violation. 
More specifically, let $\mathbf{a}^T_k$ be the $k$th row of matrix $\mathbf{A}$ in \eqref{eq:ip}, 
the constraint violation (CV) loss is defined as
$\mathcal{L}_{\text{CV}} \triangleq \frac{1}{m}\sum_{k=1}^m \max(\mathbf{a}_k^T\hat{\xx}-b_k, 0)$, 
where $m$ is the number of constraints.
The total loss for training diffusion and decoder therefore consists of the three parts:
\begin{equation}
\label{loss:total}
\mathcal{L} = \mathcal{L}_{\text{MSE}}+\mathcal{L}_{\text{CE}}+\lambda \mathcal{L}_{\text{CV}},
\end{equation}
where $\lambda$ is a hyper-parameter to regulate the penalty. 
The full training procedure is given in Algorithm \ref{alg:training} 
and the training details can be found in Appendix~\ref{app:training}.
\begin{algorithm}[ht]
\caption{Training diffusion and solution decoder}
\label{alg:training}
\begin{algorithmic}[1] 
\ENSURE IP instance embedding $\bzi$ from CISP, solution embedding $\bzx$ from CISP
\REQUIRE diffusion model $\ff_{\theta}$, and solution decoder $\dd_{\phi}$.
    \REPEAT
    \STATE $t\sim \text{Uniform}(\{1,...,T\})$ 
        \STATE $\bepsilon \sim \mathcal{N}(\mathbf{0},\bI)$
        \STATE $\hatbzx \leftarrow \ff_{\theta}\left(\sqrt{\Bar{\alpha}_t}\bzx+\sqrt{1-\Bar{\alpha}_t}\bepsilon, \bzi,t\right)$ {\color{gray}\small \textit{// predicted $\hatbzx$}}
        \STATE $\hat{\xx} \leftarrow \dd_{\phi}(\hatbzx, \bzi)$ {\color{gray}\small \textit{// reconstructed solution $\hat{\xx}$}}
        \STATE Take gradient descent step to minimize total loss in \eqref{loss:total}
        \UNTIL Reaches a fixed number epochs or satisfies an early stopping criteria
\end{algorithmic}
\end{algorithm}

\section{IP Guided Sampling}
\label{sec:guided-sampling}
Once the models have been trained, 
we can then sample variable assignments by running the sampling algorithm of DDPM or DDIM from a random Gaussian noise $\bzxT\sim \mathcal{N}(\mathbf{0}, \bI)$. 
Interestingly we find that, without suitable guidance, 
diffusion model is prone to generate inaccurate distributions, e.g. violating constraints for a given IP instance as shown in section~\ref{sec:toy_example}. 
We thus consider the constraints information $(\mathbf{A}, \mathbf{b})$ and objective coefficient $\mathbf{c}$ during sampling. 
We present the IP guided diffusion sampling for both DDPM and DDIM, 
of which the latter is faster and better in terms of the quality and feasibility, as shown in Section~\ref{sec:experiments}.

\subsection{IP Guided Diffusion Sampling}
\label{sec:IP-ddpm}
Consider a conditional diffusion model $p_{\theta}(\bzxt|\bzx^{(t+1)}, \bzi)$, we first introduce \textit{constraint guidance} by designing each transition probability as 
\begin{equation}
\label{eq:ddpm_transition}
    p_{\theta, \phi}(\bzxt|\bzx^{(t+1)}, \bzi, \mathbf{A}, \mathbf{b}) = Zp_{\theta}(\bzxt|\bzx^{(t+1)}, \bzi) e^{-sc_{\phi}(\bzx^{(t)},\bzi, \mathbf{A}, \mathbf{b})},
\end{equation}
where $s$ is the gradient scale, $Z$ is a normalizing constant and $c_{\phi}(\bzxt, \bzi, \mathbf{A},\mathbf{b}) = \sum_{k=1}^m \max (\mathbf{a}_k^T\dd_{\phi}(\bzxt, \bzi)-b_k,0)$ measures the violation of constraints.
Let $\bmu$ and $\bSigma$ be the mean and variance of the Gaussian distribution representing $p_{\theta, \phi}(\bzxt|\bzx^{(t+1)}, \bzi)$. Then,
$\log p_{\theta}(\bzxt|\bzx^{(t+1)}, \bzi) = -\frac{1}{2}(\bzxt-\bmu)^T\bSigma^{-1}(\bzxt-\bmu) + C$, 
where $C$ is a constant. 
Consider the Taylor expansion for $c_{\phi}$ at $\bzxt = \bmu$,
\begin{align*}
    &c_{\phi}(\bzxt,\bzi, \mathbf{A}, \mathbf{b}) 
    \\
    &\approx c_{\phi}(\bmu, \bzi, \mathbf{A}, \mathbf{b})+(\bzxt-\bmu) \nabla_{\bzxt} c_{\phi}(\bzxt, \bzi, \mathbf{A}, \mathbf{b})|_{\bzxt=\bmu}\\
    &=  (\bzxt-\bmu)\mathbf{g} +C_1,
\end{align*}
where $\mathbf{g} = \nabla_{\bzxt} c_{\phi}(\bzxt, \bzi, \mathbf{A}, \mathbf{b})|_{\bzxt=\bmu}$ and $C_1$ is a constant.
Similar to Classifier Guidance~\citep{dhariwal2021diffusion}, 
we assume that $c_{\phi}(\bzxt, \bzi, \mathbf{A}, \mathbf{b})$ has low curvature compared to $\bSigma^{-1}$ and thus have the following,
\begin{align*}
    & \log (p_{\theta, \phi}(\bzxt|\bzx^{(t+1)}, \bzi, \mathbf{A}, \mathbf{b})) \\
    &\approx -\frac{1}{2}(\bzxt-\bmu)^T\bSigma^{-1}(\bzxt-\bmu)- s(\bzxt-\bmu)\mathbf{g} + sC_1 + C_2 \nonumber \\
    &= -\frac{1}{2}(\bzxt-\bmu+s\bSigma \mathbf{g})^T\bSigma^{-1}(\bzxt-\bmu+s\bSigma \mathbf{g})+C_3,
\end{align*}
where $C_3$ is a constant and can be safely ignored. Therefore, the denoising transition $p_{\theta, \phi}(\bzxt|\bzx^{(t+1)}, \bzi, \mathbf{A}, \mathbf{b})$ can be approximated by a Gaussian distribution with a mean shifted by $-s\bSigma \mathbf{g}$.
We can further inject the objective guidance to transition probability for acquiring high-quality solutions,
\begin{equation}
\begin{aligned}
    \label{eq:ddpm sampling}
    &p_{\theta, \phi}(\bzxt|\bzx^{(t+1)}, \bzi, \mathbf{A}, \mathbf{b}, \mathbf{c}) \\
    &= Zp_{\theta}(\bzx^{(t)}|\bzx^{(t+1)}, \bzi) e^{-s\left((1-\gamma)c_{\phi}(\bzx^{(t)},\bzi, \mathbf{A}, \mathbf{b}) + \gamma o_{\phi}(\bzx^{(t)},\bzi, \mathbf{c})\right)},
\end{aligned}
\end{equation}
where $o_{\phi}(\bzxt,\bzi, \mathbf{c}) = \mathbf{c}^T\dd_{\phi}(\bzxt, \bzi)$, $\mathbf{c}$ is the coefficient of objective from \eqref{eq:ip}, 
and $\gamma\in [0,1]$ is the leverage factor for balancing constraint and objective. The corresponding sampling method is called \emph{IP Guided Diffusion Sampling}, as presented in Algorithm~\ref{alg:ddpm}.

\begin{algorithm}[ht]
\caption{IP Guided Diffusion Sampling}
\label{alg:ddpm}
\begin{algorithmic}[1] 
\ENSURE gradient scale $s$, leverage factor $\gamma$, constraint information $(\mathbf{A},\mathbf{b})$ and objective coefficient $\mathbf{c}$ 
\REQUIRE diffusion model $\ff_{\theta}$ and solution decoder $\dd_{\theta}$.
\STATE $\bzx^{(T)}\sim \mathcal{N}(\mathbf{0},\mathbf{I})$
\FOR{ $t$ from $T$ to $1$}
\STATE $\bmu \leftarrow \frac{\sqrt{\alpha_t}(1-\Bar{\alpha}_{t-1})}{1-\Bar{\alpha}_t}\bzxt + \frac{\sqrt{\Bar{\alpha}_{t-1}}\beta_t}{1-\Bar{\alpha}_t}\ff_{\theta}(\bzxt,\bzi,t)$
\STATE $\bSigma \leftarrow  \bSigma_{\theta}(\bzxpt|\bzxt, \bzi)$
\STATE $\bmu \leftarrow \bmu - s\bSigma\nabla_{\bzxt} \left((1-\gamma)c_{\phi}(\bzxt, \bzi, \mathbf{A}, \mathbf{b}) + \gamma o_{\phi}(\bzxt,\bzi, \mathbf{c})\right)$
\STATE $\bzxpt \sim \mathcal{N}(\bmu, \bSigma)$
\ENDFOR
\RETURN $\dd_{\phi}(\bzxzero, \bzi)$.
\end{algorithmic}
\end{algorithm}

\subsection{Non-Markovian IP Guided Sampling}
\label{sec:IP-ddim}
For the non-Markovian sampling scheme as used in DDIM, 
the method for conditional sampling is no longer invalid. 
To guide the sampling process, 
existing studies~\citep{ho2022classifier, dhariwal2021diffusion} used the score-based conditioning trick from~\citet{song2020score} to construct a new epsilon prediction. 
We found that this trick turns out ineffective in our problem setting, 
possibly due to unique difficulties of finding feasible solutions in IP instances. 

Instead, we find that adding a direction that guides the predicted solution to constraint region in each step of reverse process helps generate more reasonable solutions. 
Specifically, we first generate the predicted noise according to $\hatbzx^{(0)} = f_{\theta}(\bzxt, \bzi, t)$, 
$\bepsilon_{\theta}^{(t)} = \frac{\bzxt-\sqrt{\Bar{\alpha}_t}\hatbzx^{(0)} }{\sqrt{1-\Bar{\alpha}_t}}$. 
According to~\eqref{eq:ddim-reverse}, the transition equation of DDIM for $\bzxpt$ from a sample $\bzxt$ can be written as 
\begin{equation}
\label{eq:ddim-transition}
    \bzxpt = \sqrt{\Bar{\alpha}_t} \ff_{\theta}(\bzxt, \bzi, t) + \sqrt{1-\Bar{\alpha}_{t-1}-\sigma_t^2} \bepsilon_{\theta}^{(t)}+\sigma_t\bepsilon_t.
\end{equation}
where the first term is the prediction of $\bzxzero$ and the second term is the direction pointing to $\bzxt$. 
To consider constraints in~\eqref{eq:ip}, we modify $\bepsilon_{\theta}^{(t)}$ by adding the direction of minimizing sum of constraint violation, that is
\begin{equation}
    \label{eq:ddim-cg}
    \hat{\bepsilon}_{\theta}^{(t)} = \bepsilon_{\theta}^{(t)}- s\nabla_{\bzxt} c_{\phi}(\bzxt, \bzi, \mathbf{A}, \mathbf{b}),
\end{equation}
where $s$ is gradient scale. 
By replacing $\bepsilon_{\theta}^{(t)}$ in \eqref{eq:ddim-transition} with $\hat{\bepsilon}_{\theta}^{(t)}$, 
we obtain \emph{Non-Makovian Constraint Guided Sampling}, 
which guides the solution generated in each transition to approach constraint region. 
Equivalently, it is to perform a gradient descent step with a step-size shrinking to zero as $t \to 0$ when $\sigma_t\to 0$. 
To further consider the objective function together with the constraint, we can update $\bepsilon_{\theta}^{(t)}$ as follows:
\begin{equation}
    \label{eq:ddim-ipg}
    \hat{\bepsilon}_{\theta}^{(t)} = \bepsilon_{\theta}^{(t)}- s\nabla_{\bzxt}\bigg((1-\gamma) c_{\phi}(\bzxt, \bzi, \mathbf{A}, \mathbf{b}) + \gamma o_{\phi}(\bzxt, \bzi, \mathbf{c})\bigg),
\end{equation}
where $\gamma\in [0,1]$ is the leverage factor for balancing constraint and objective. 
This method is called \textit{Non-Markovian IP Guided Diffusion Sampling}, as presented in Algorithm~\ref{alg:ddim} .
\begin{algorithm}[ht]
\caption{Non-Markovian IP Guided Diffusion Sampling }
\label{alg:ddim}
\begin{algorithmic}[1] 
\ENSURE gradient scale $s$, leverage factor $\gamma$, constraint information $(\mathbf{A},\mathbf{b})$ and objective coefficient $\mathbf{c}$
\REQUIRE diffusion model $\ff_{\theta}$, solution decoder $\dd_{\theta}$
    \STATE $\bzxT \sim \mathcal{N}(\mathbf{0},\bI)$
\FOR{ $t$ from $T$ to $1$}
    \STATE $\bepsilon_{\theta}^{(t)} \leftarrow (\bzxt-\sqrt{\Bar{\alpha}_t}\ff_{\theta}(\bzxt, \bzi, t))/\sqrt{1-\Bar{\alpha}_t}$.
\STATE $\hat{\bepsilon}_{\theta}^{(t)} \leftarrow \bepsilon_{\theta}^{(t)}- s\nabla_{\bzxt} \bigg((1-\gamma)c_{\phi}(\bzxt, \bzi, \mathbf{A}, \mathbf{b}) + \gamma o_{\phi}(\bzxt, \bzi, \mathbf{c})\bigg)$
    \STATE $\bzxpt \leftarrow \sqrt{\Bar{\alpha}_t} \ff_{\theta}(\bzxt, \bzi, t) + \sqrt{1-\Bar{\alpha}_{t-1}-\sigma_t^2} \hat{\bepsilon}_{\theta}^{(t)}+\sigma_t\bepsilon_t$
    \ENDFOR
\RETURN $\dd_{\phi}(\bzxzero, \bzi)$.
\end{algorithmic}
\end{algorithm}

\begin{figure*}[ht]
    \centering
    \includegraphics[width=0.7\linewidth]{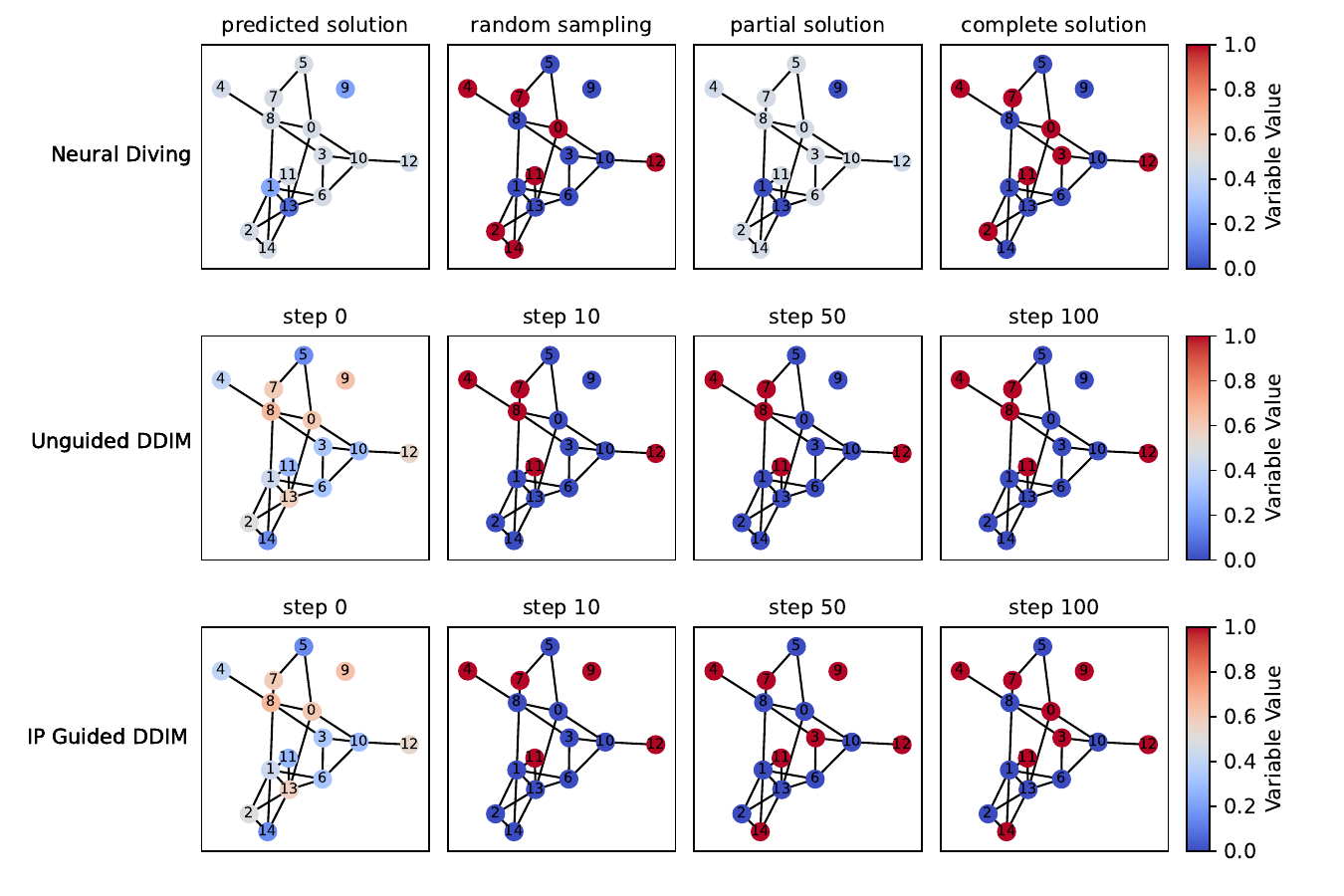}
    \caption{The sampling results from different methods. For Neural Diving, we present the predicted solution from GCN, random sampling according to the predicted solution, the partial solution obtained via SelectiveNet (only node 1, 9 and 13 are assigned to 0), and the completing result by calling CompleteSol heuristic. For DDIM and IP Guided DDIM, we present the results from different time steps (transformed to solution space by a decoder) during sampling.}
    \label{fig:toy_example}
\end{figure*}

\section{Experiments}
\label{sec:experiments}

This section empirically investigates the effectiveness of our method in solving IP instances. 
The efficacy is evaluated with two metrics: \emph{feasible ratio} and \emph{objective value}. 
The feasible ratio measures the proportion of feasible solutions among all sampled solutions and the objective value obtained from the generated feasible solutions measures the solution quality. Additionally, we also present gap between feasible solution $x$ and optimal solution $x^*$ via 
$\frac{|\mathbf{c}^T(\mathbf{x}-\mathbf{x}^*)|}{\max(|\mathbf{c}^T\mathbf{x}|, |\mathbf{c}^T\mathbf{x}^*|)}$ for comparing the performance more clearly.

We evaluate our methods on four IP datasets generated by the Ecole library~\citep{prouvost2020ecole}: 
\begin{itemize}[noitemsep,topsep=0pt,parsep=0pt,partopsep=0pt,leftmargin=*]
    \item \emph{Set Cover} (SC) is to find the least number of subsets that cover a given universal set. 
    \item \emph{Capacitated Facility Location} (CF) is to locate a number of facilities to serve the sites with a given demand and the aim is minimize the total cost. 
    \item \emph{Combinatorial Auction} (CA) is to help bidders place unrestricted bids for bundles of goods and the aim is to maximize the revenue.
    \item \emph{Independent Set} (IS) is to find the maximum subset of nodes of an undirected graph such that no pair of nodes are connected.
\end{itemize}
 We compare the performance of our approach with state-of-the-art methods including Neural Diving (ND) \citep{nair2020solving}, the Predict-and-Search algorithm (PS) \citep{han2023a}, and heuristic solutions from SCIP 8.0.1 \citep{bestuzheva2021scip} and Gurobi 9.5.2 \citep{gurobi2021gurobi}. Detailed descriptions of the datasets and baseline methods are provided in Appendix~\ref{app:data_baselines}. Moreover, we incorporate the best objective values, which serve as the ground truth, acquired by executing Gurobi on each instance for a duration of 100 seconds. To ensure clarity, we use \emph{IP Guided DDPM} to denote the (Markovian) IP Guided Diffusion sampling in Section \ref{sec:IP-ddpm}, and \emph{IP Guided DDIM} to represent the Non-Markovian IP Guided Diffusion sampling in Section \ref{sec:IP-ddim}.  

In the following, we first illustrate the guided diffusion sampling and emphasize its distinctions to Neural Diving and vanilla generation process in diffusion models in Section~\ref{sec:toy_example}. 
Further, we evaluate the feasibility and quality of solutions generated by IP guided DDIM across all four datasets in Section~\ref{sec:compare}. Furthermore, we conduct an ablation study to investigate the impact of different guided approaches and contrastive learning in Section~\ref{sec:ablation_study}.
In Section~\ref{sec:generalization}, we demonstrate the scalability of our approach by applying it to larger instances and present the outcomes of qualitative analysis of solutions.
All experiments are performed in a workstation with two Intel(R) Xeon(R) Platinum 8163 CPU @ 2.50GHz, 176GB ram and two Nvidia V100 GPUs. 
We also provide the total training and inference time in Appendix~\ref{sec:train_inf_time}. The detailed hyper-parameters for IP guided sampling can be found in Appendix~\ref{app:sampling-parameter}. 
The codes can be found in: \url{https://github.com/agent-lab/diffusion-integer-programming}.

\begin{figure*}[ht]
\centering 
    \includegraphics[width=0.7\textwidth]{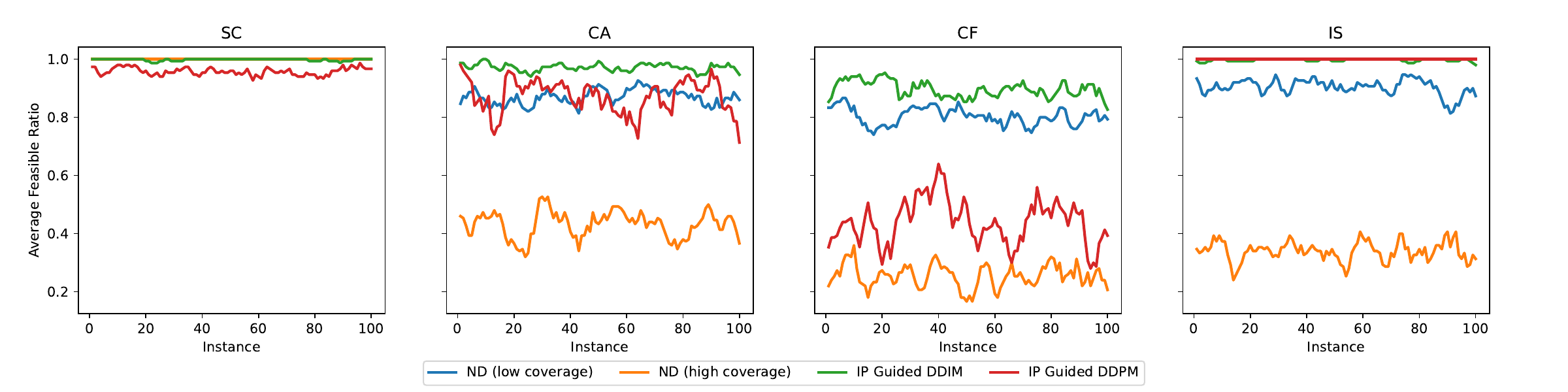}
\caption{The feasible ratio in 100 instances, with each instance sampled 30 complete or partial solutions. For diffusion model, we measure the feasible ratio of  complete solutions. Two versions of Neural Diving are trained with distinct coverage thresholds, referred to as ND (low coverage) and ND (high coverage). The feasibility ratio is evaluated only for partial solutions from Neural Diving, as the complete solutions from this method yield a 0\% feasibility ratio.}
\label{fig:feasibility} 
\end{figure*}

\begin{table*}[ht]
\caption{The average objective value (obj.), gap and feasible ratio (fea.) for 100 instances on 4 datasets. Optimal, heuristic, complete and partial denote optimal solutions from Gurobi, heuristic solutions from solvers, complete solutions from models, and partial solutions from models. CompleteSol and search indicates the CompleteSol heuristic and predict-and-search algorithm for completing the partial solutions, respectively.}
\resizebox{1\textwidth}{!}{
\begin{tabular}{c|c|ccc|ccc|ccc|ccc}
\toprule
\multirow{2}{*}{Algorithm} & \multirow{2}{*}{Solutions' Type} & \multicolumn{3}{c|}{SC (min)}                     & \multicolumn{3}{c|}{CF (min)}                       & \multicolumn{3}{c|}{CA (max)}                       & \multicolumn{3}{c}{IS (max)}                      \\
                           &                       & obj.           & gap             & fea.           & obj.             & gap             & fea.           & obj.             & gap             & fea.           & obj.           & gap             & fea.           \\ \midrule
Gurobi (100s)                     & Optimal                 & 168.3          & 0.0\%           & 100\%          & 11405.5          & 0.0\%           & 100\%          & 36102.6          & 0.0\%           & 100\%          & 685.3          & 0.0\%           & 100\%          \\ \midrule
SCIP                       & Heuristic                  & 1967.0         & 91.4\%          & \textbf{100\%} & 84748.4          & 86.5\%          & \textbf{100\%} & 28007.4          & 22.4\%          & \textbf{100\%} & 447.8          & 34.7\%          & \textbf{100\%} \\
Gurobi                     & Heuristic                  & \textbf{522.4} & \textbf{67.8\%} & \textbf{100\%} & 50397.3          & 77.4\%          & \textbf{100\%} & \textbf{30052.0} & 16.8\%          & \textbf{100\%} & 415.5          & 39.4\%          & \textbf{100\%} \\
Neural Diving              & Complete                    & -              & -               & 0.0\%          & -                & -               & 0.0\%          & -                & -               & 0.0\%          & -              & -               & 0.0\%          \\
IP Guided DDPM (Ours)      & Complete                    & 577.9          & 70.8\%          & 95.7\%         & 58488.1          & 80.5\%          & 44.0\%         & 800.3            & 97.8\%          & 87.3\%         & 129.9          & 81.0\%          & \textbf{100\%} \\
IP Guided DDIM (Ours)      & Complete                    & 533.5          & 68.5\%          & 99.8\%         & \textbf{25119.2} & \textbf{54.6\%} & 89.7\%         & 26916.9          & 25.4\%          & 97.1\%         & \textbf{455.6} & \textbf{33.5\%} & 99.7\%         \\ \midrule
Neural Diving              & Partial + CompleteSol            & 849.0          & 80.2\%          & \textbf{100\%} & 14259.8          & 20.0\%          & 81.3\%         & 30143.6          & 16.5\%          & 87.0\%         & 484.1          & 29.4\%          & 90.4\%         \\
PS                         & Partial + Search            & 593.7          & 71.7\%          & \textbf{100\%} & 32119.8          & 64.5\%          & \textbf{100\%} & 31159.5          & 13.7\%          & \textbf{100\%} & 587.9          & 14.2\%          & \textbf{100\%} \\
IP Guided DDIM (Ours)      & Partial + CompleteSol            & \textbf{255.5} & \textbf{34.1\%} & \textbf{100\%} & \textbf{14224.1} & \textbf{19.8\%} & \textbf{100\%} & \textbf{32491.1} & \textbf{10.0\%} & 99.7\%         & \textbf{639.4} & \textbf{6.7\%}  & \textbf{100\%} \\ \bottomrule
\end{tabular}
}
 \label{tab:compare}
\end{table*}

\subsection{Illustrative experiments}
\label{sec:toy_example}

The maximal independent set problem involves finding the largest subset of nodes in an undirected graph where no two nodes are connected. 
We focus on an illustrative example: a graph consisting of 15 nodes and 22 edges.  
This graph can be transformed into an integer programming (IP) instance with 15 variables and 22 constraints. The graph is depicted in Figure \ref{fig:toy_example}, and we present the results from Neural Diving, and from the different time steps of unguided DDIM and IP Guided DDIM. 
Among the three algorithms, Neural Diving fixes only three node values and uses the Completesol heuristic to find an independent set containing seven nodes. 
However, the random sampling solution based on predicted probability from Neural Diving is infeasible. 
Unguided DDIM is unable to find a feasible solution (there is an edge between node 7 and node 8).
In contrast, IP Guided DDIM is able to fetch the optimal solution by finding an independent set containing eight nodes during the sampling process, where no two nodes are connected. 
Notably, the quality of the solution improves as the sampling process progresses. 
The independent set contains 5 nodes at step 10, 7 nodes at step 50, and finally 8 nodes (the optimal solution) at step 100. 
These indicate that IP Guided DDIM outperforms Neural Diving and Unguided DDIM in finding the optimal solution for this illustrative example.

\subsection{Performance Evaluation}
\label{sec:compare}

\begin{table*}[ht]
\centering
\caption{Ablation study for 100 instances on 4 datasets with different guidances.}
\resizebox{0.7\linewidth}{!}{
\begin{tabular}{lcccccccccc}
\hline
\multicolumn{1}{c}{} &
  \multicolumn{2}{c}{\begin{tabular}[c]{@{}c@{}}\small Unguided \\ \small DDIM \end{tabular}} &
  \multicolumn{2}{c}{\begin{tabular}[c]{@{}c@{}}\small Constraint Guided \\ \small DDIM\end{tabular}} &
  \multicolumn{2}{c}{\begin{tabular}[c]{@{}c@{}}\small Objective Guided \\ \small DDIM \end{tabular}} &
  \multicolumn{2}{c}{\begin{tabular}[c]{@{}c@{}}\small IP Guided DDIM \\ \small w/o CISP \end{tabular}} &
  \multicolumn{2}{c}{\begin{tabular}[c]{@{}c@{}} \small IP Guided \\ \small DDIM \end{tabular}} \\ \cline{2-11} 
\multicolumn{1}{c}{\multirow{-2}{*}{\small dataset}} 
& \small obj. & \small  fea. & \small obj. & \small fea. & \small obj.    & \small fea. & \small  obj.    & \small fea. & \small  obj.    & \small fea.\\ \hline
\small SC (min)   & -   & \small  0.0\%  &  \small 63046.9 &  \small \textbf{99.8\%} & -    & \small  0.0\%  & \small 763.4 & \small 99.8\% &  \small \textbf{533.5}   &  \small \textbf{99.8\%} \\
\small CF (min)   & -    &  \small 0.0\%  &  \small 53311.2 &  \small 74.1 \%& -    &  \small 0.0\%  &  \small 31319.8 & \small 41.7\% & \small \textbf{25119.2} &  \small \textbf{89.7\%} \\
\small CA (max)   & -    & \small  0.0\%  &  \small 5157.2  &  \small \textbf{99.7\%} & -    & \small  0.0\%  & \small 23383.3 & \small 57.7\% & \small \textbf{26916.9} &  \small 97.1 \%\\
\small IS (max)  & -    &  \small 0.0\%  &  \small 386.5   &  \small \textbf{100 \%} & -    &  \small 0.0\%  &   \small \textbf{479.1} & \small 68.9\% & \small 455.6  &  \small 99.7\% \\ \hline
\end{tabular}}
\label{tab:ablation}
\end{table*}

\begin{table*}
\centering
 \caption{The average objective value (obj.), gap and the feasible ratio (fea.) for 100 instances in 3 different size SC datasets. Optimal, heuristic, complete and partial denote optimal solutions from Gurobi, heuristic solutions from solvers, complete solutions from models, and partial solutions from models. CompleteSol and search indicates CompleteSol heuristic and predict-and-search algorithm for completing the partial solutions, respectively.}
\resizebox{0.8\textwidth}{!}{
\begin{tabular}{c|c|ccc|ccc|ccc}
\toprule
\multirow{2}{*}{Algorithm} & \multirow{2}{*}{Solutions' Type} & \multicolumn{3}{c|}{SC (2000)}                    & \multicolumn{3}{c|}{SC (3000)}                    & \multicolumn{3}{c}{SC (4000)}                     \\
                           &                       & obj.           & gap             & fea.           & obj.           & gap             & fea.           & obj.           & gap             & fea.           \\ \midrule
Gurobi (100s)                     & Optimal                 & 168.3          & 0.0\%           & 100\%          & 140.4          & 0.0\%           & 100\%          & 126.9          & 0.0\%           & 100\%          \\ \midrule
SCIP                       & Heuristic                  & 1977.0         & 91.4\%          & 100\%          & 2236.2         & 93.7\%          & \textbf{100\%} & 2386.0         & 94.7\%          & \textbf{100\%} \\
Gurobi                     & Heuristic                  & 522.4          & 67.8\%          & 100\%          & 718.8          & 80.5\%          & \textbf{100\%} & 1454.5         & 91.3\%          & \textbf{100\%} \\
Neural Diving              & Complete                    & -              & -               & 0.0\%          & -              & -               & 0.0\%          & -              & -               & 0.0\%          \\
IP Guided DDPM (Ours)      & Complete                    & 594.7          & 70.8\%          & 96.5\%         & \textbf{451.8} & \textbf{68.9\%} & 83.7\%         & \textbf{440.7} & \textbf{71.2\%} & 77.9\%         \\
IP Guided DDIM (Ours)      & Complete                    & \textbf{533.5} & \textbf{68.5\%} & 99.8\%         & 486.8          & 71.2\%          & 99.9\%         & 464.9          & 72.7\%          & \textbf{100\%} \\ \midrule
Neural Diving              & Partial + CompleteSol            & 849.0          & 80.2\%          & \textbf{100\%} & 1145.8         & 87.7\%          & \textbf{100\%} & 1465.6         & 91.3\%          & \textbf{100\%} \\
PS                         & Partial + Search            & 593.7          & 71.7\%          & \textbf{100\%} & 737.0          & 80.9\%          & \textbf{100\%} & 994.9          & 87.2\%          & \textbf{100\%} \\
IP Guided DDIM (Ours)      & Partial + CompleteSol            & \textbf{255.5} & \textbf{34.1\%} & \textbf{100\%} & \textbf{217.4} & \textbf{35.4\%} & \textbf{100\%} & \textbf{195.9} & \textbf{35.2\%} & \textbf{100\%} \\ \bottomrule
\end{tabular}}
 \label{table:generalization}
\end{table*}

In this section, we evaluate the performance of different methods by comparing their average feasible ratios, gaps and average objective values across four datasets mentioned earlier. Each dataset contains 100 instances. For each instance, we sample 30 solutions and calculate the corresponding metrics, which allows us to assess the performance of each method in terms of both solution feasibility and objective value across the different datasets.

We first compare the feasibility ratio of the \emph{complete} solutions generated by IP Guided Diffusion with the \emph{partial} solutions generated by Neural Diving, due to the inability of Neural Diving to produce complete solutions. In Neural Diving, the expected proportion of variables assigned by the model is controlled by the hyper-parameter "Coverage" $C$. A higher $C$ results in more variables being assigned by the model, which generally leads to a lower feasibility ratio. Therefore, a carefully prescribed coverage $C$ is crucial in Neural Diving to ensure the feasibility of partial solutions.
In this experiment, we trained two variants of Neural Diving with different coverage thresholds for each dataset. 
For the SC, CA, and IS datasets, the coverage are set to 0.2 (low coverage) and 0.3 (high coverage) respectively. However, for the CF dataset, the coverage thresholds are set to 0.1 and 0.2.
Figure~\ref{fig:feasibility} presents the average feasible ratio for 100 instances. In this comparison, the feasible ratio of complete solutions from IP guided DDIM outperforms the solutions of Neural Diving in almost all instances.

To more comprehensively evaluate the performance of the diffusion model, we compare it against different baselines in Table~\ref{tab:compare}. For SCIP, we adopt the first solution obtained through non-trivial heuristic algorithms during the solving phase. For Gurobi, we use the best heuristic solution for each instance as a benchmark. Additionally, we include the optimal objective values (obtained through running Gurobi for 100 seconds on each instance) as ground-truth. For Neural Diving, complete solutions are frequently infeasible. As such, we employ a low-coverage model that prioritizes the feasibility of partial solutions. Subsequently, a CompleteSol heuristic is utilized to finalize these partial solutions. In the case of the Predict and Search algorithm (PS), we construct a trust region using partial solutions with the same proportion of assigned variables as Neural Diving and use Gurobi as the Solver to search the best heuristic solutions found as a benchmark. To ensure a fair comparison, we report not only the quality of the complete solutions generated by our methods (IP Guided DDIM/DDPM), but also the quality of partial solutions. For this, we randomly select the same proportion of variables from DDIM's complete solutions as the ratio from Neural Diving and PS as our partial solutions. We then use the CompleteSol heuristic to finalize these partial solutions.

The results are presented in Table~\ref{tab:compare}. It is evident that the complete solutions produced using the IP Guided DDIM method have a feasible ratio of at least 89.7\%, and their objective values are comparable to the best heuristic solutions from Gurobi. In contrast, the complete solutions obtained solely through Neural Diving are consistently infeasible. Moreover, the integration of partial solutions from IP Guided DDIM with the CompleteSol heuristic outstrips that of all methods in terms of objective values. This improvement is demonstrated by a 3.7 to 33.7\% reduction in the gap to optimal values, while the feasibility ratio for all datasets approaches 100\%.

\subsection{Ablation Study}
\label{sec:ablation_study}

We ablate on unguided DDIM (with $s = 0$), constraint guided DDIM (with $\gamma = 1$), objective guided DDIM (with $\gamma = 0$) models and IP guided DDIM on four datasets. 
We also include an experiment where we train IP and solution embeddings directly via algorithm \ref{alg:training} without CISP, in order to assess the advantages of contrastive learning, i.e. IP Guided DDIP w/o CISP in Table~\ref{tab:ablation}.   
The results are presented in Table \ref{tab:ablation}. 
Evidently, the constraint guidance is crucial in generating feasible solutions, 
and the objective guidance further enhances the quality of solutions. 
Moreover, the experiments demonstrate that CISP plays a crucial role in ensuring that the solutions produced by our methods are more feasible.
Therefore, combining both constraint and objective guidance achieves good quality solutions with high probability. 

\begin{figure}
\centering
\subfigure[SC instance (minimization).]{\includegraphics[width=0.35\textwidth]{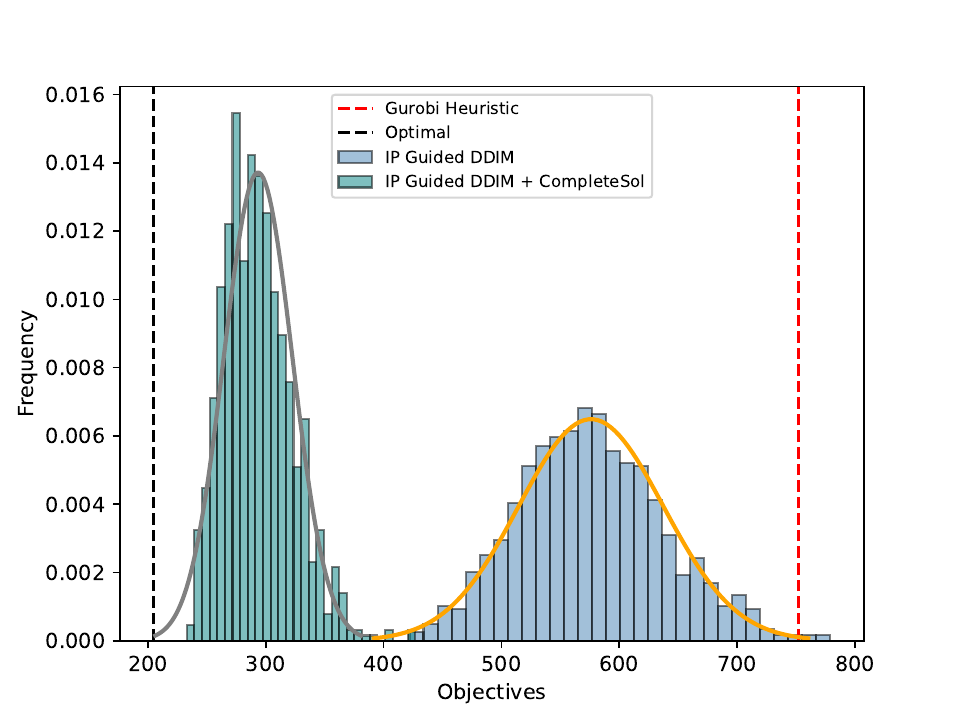}}
\subfigure[IS instance (maximization).]{\includegraphics[width=0.35\textwidth]{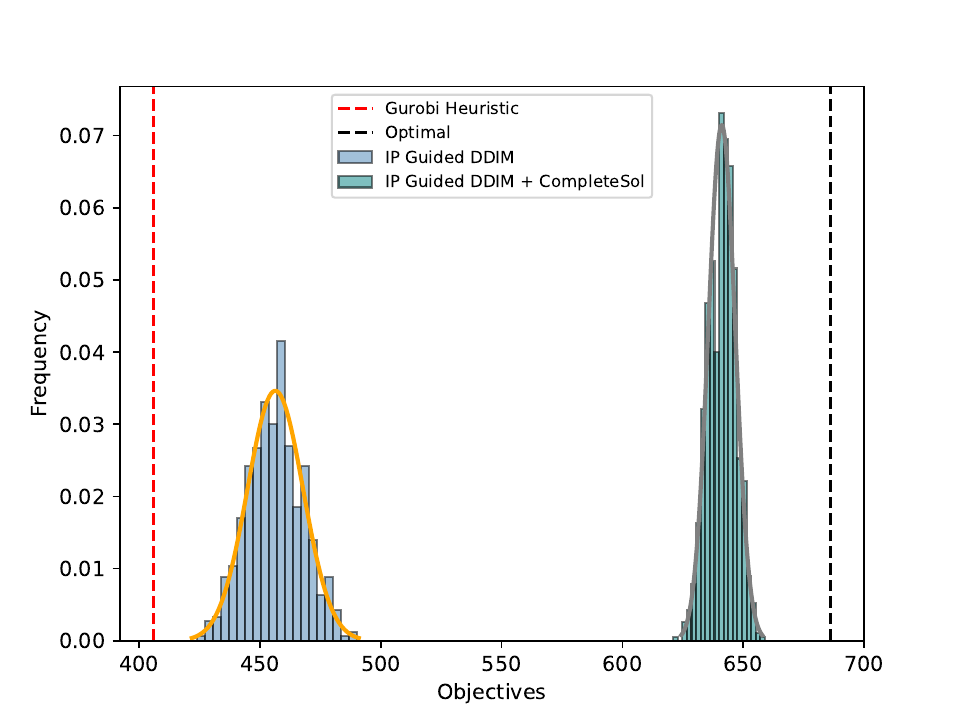}}
\caption{The objective distribution of 1000 solutions sampled from a single instance.}
\label{fig:distribution_of_solutions}
\end{figure}

\subsection{Scalability test and qualitative analysis}
\label{sec:generalization}

Practitioners often aim to apply the models learned to solve problems of larger scales than the ones used for data collection and training. 
To estimate how well a model can generalize to bigger instances, we evaluate its performance on datasets of varying sizes from the Set Cover problem (minimization problem). We utilize three size categories: 
\begin{itemize}
    \item SC (2000): 2000 variables, 1000 constraints
    \item SC (3000): 3000 variables, 1500 constraints
    \item SC (4000): 4000 variables, 2000 constraints
\end{itemize}
It is worth noting that all models were trained using the SC (2000) dataset. The results in table \ref{table:generalization} demonstrate that IP guided DDPM consistently performs well across all three different-sized datasets, indicating that our framework possesses strong generalization capabilities.

Diffusion models are generative models that capture the distribution of a dataset. In this experiment, we focus on the distribution of solutions generated by our methods. We take a single instance from the SC dataset and IS dataset and use the IP Guided DDIM algorithm to generate 1000 complete feasible solutions. We also randomly sample 20\% of variable values from each solution and use the CompleteSol heuristic to fill in the remaining variables. We analyze the distribution of objective values for these complete solutions and compare them with the optimal objective values and the best heuristic solutions found by the Gurobi heuristic.
The results are shown in Figure \ref{fig:distribution_of_solutions}. 
Clearly, the complete solutions directly generated by the IP Guided DDIM algorithm are superior to the best solutions from the Gurobi heuristic. 
Furthermore, the objective values from the partial solutions completed by the CompleteSol heuristic are closer to the optimal value.

\section{Related Work}
\label{sec:related_works}
We start the related work section with a summary of deep learning techniques used in the construction of feasible solutions for Integer Programming (IP) problems. 
\citet{gasse2019exact} propose a method that combines a bipartite graph with a Graph Convolutional Network (GCN) to extract representations of Integer Programming (IP) instances. 
Although this approach is primarily employed to learn the branching policy in the branch and bound algorithm, 
it is worth noting that this modeling method can also be utilized for the prediction of solutions for IP instances.
However, the solutions produced by GCN directly are often infeasible or sub-optimal.
To address this, Neural Diving~\citep{nair2020solving} leverages the SelectiveNet~\citep{geifman2019selectivenet}
to assign values to only a subset of the variables based on a coverage threshold, 
with the rest of the variables being completed via an IP Solver.
To further improve the feasibility of generated solutions, 
\citet{han2023a} proposes a predict-search framework that combines the predictions from GNN model with trust region method. 
However, this method still relies on IP solver to solve a modified instance (adding neighborhood constraints to origin instance) in order to get complete solutions. 
Similar to \citet{nair2020solving}, \citet{khalil2022mip} integrate GNN into integer programming solvers and apply it to construct a partial solutions through a prescribed rounding threshold, 
which is then completed using SCIP. 
In contrast, our method aims to learn the latent structure of IP instances by diffusion models, 
and obtains complete feasible solutions through guided diffusion sampling, without any reliance on the IP solver.

Another set of related works to our paper is diffusion models~\citep{sohn2015learning, ho2020denoising}. 
As the latest state-of-the-art family of deep generative models, 
diffusion models have demonstrated their ability to enhance performance across various generative tasks~\citep{yang2022diffusion}. 
In this paper, we focus our discussion specifically on conditional diffusion models. 
Unlike unconditional generation, conditional generation emphasizes application-level contents as a condition to control the generated results based on predefined intentions.
To enable this conditioning, \citet{dhariwal2021diffusion} introduce the concept of classifier guidance, 
which enhances sample quality by conditioning the generative process on an additional trained classifier. 
In the same vein, \citet{ho2022classifier} propose a joint training strategy for both conditional and unconditional diffusion models, 
i.e., classifier-free guidance. 
This approach combines the resulting conditional and unconditional scores to achieve a balance between sample quality and diversity. 
This idea has also found its effectiveness in Topology Optimization~\citep{maze2023diffusion}. 
We thus take a similar derivation from the classifier guidance and devise IP-guided diffusion sampling by incorporating the objectives and constraints into the transition probability.

\section{Discussions}
\subsection{The Importance of Solution Generation}

Integer Programming solvers typically rely on branch-and-bound or branch-and-cut algorithms to tackle problems, often starting with a primal heuristic to find a good initial feasible solution. A robust starting solution generally accelerates the entire solving process. Our method demonstrates that diffusion models are capable of generating complete feasible solutions that satisfy all constraints. These generated solutions exhibit competitive performance compared to the heuristics employed in optimization tools such as Gurobi. Notably, the ability to generate complete solutions from scratch distinguishes our model from other neural-based approaches like Neural Diving \citep{nair2020solving} and the Predict-Search algorithm \citep{han2023a}. Moreover, generating feasible solutions end-to-end is especially important, as it represents a significant leap toward developing a purely neural-based method for integer programming. The paper shows that diffusion models can produce solutions that satisfy all constraints, which lays the foundation for future research in which neural networks could potentially solve such problems without relying on traditional solvers or heuristics.

\subsection{Sampling Efficiency}
The sampling time is a notable limitation of current diffusion models. Specifically, our method requires 100 iterative denoising steps for DDIM and 1000 for DDPM, resulting in longer times to generate a single solution compared to other methods. To address this issue, there have been active studies focused on accelerating diffusion sampling, such as model quantization (\citet{li2023qdiffusion}) and distillation (\citet{huang2023accelerating}). These techniques can be readily integrated into our methods. Moving forward, we will enhance our method by incorporating these acceleration techniques.

\section{Conclusion}
In this paper, we presented a comprehensive framework for generating feasible solutions for Integer Programming (IP) problems. 
We utilized the CISP approach to establish a link between IP instances and their corresponding solutions, 
allowing to obtain IP embeddings and solution embeddings. 
To effectively capture the distribution of feasible solutions, we leveraged diffusion models, which are known for their powerful learning capabilities, to learn the distribution of solution embeddings. 
We further employed a solution decoder to reconstruct the solutions from their embeddings. 
Importantly, we proposed an IP guided sampling algorithm that explicitly incorporates the objective and constraint information to generate high-quality solutions. 
The experimental results on four distinct datasets demonstrate the superiority of our approach to the state-of-the-art.

\begin{acks}
Mingfei Sun is affiliated with Centre for AI Fundamentals. He is also a member of AI Hub in Generative Models, funded by Engineering and Physical Sciences Research Council (EPSRC), part of UK Research and Innovation (UKRI),
\end{acks}

\newpage
\bibliographystyle{ACM-Reference-Format}
\balance
\bibliography{authordraft}


\begin{thebibliography}{38}


\ifx \showCODEN    \undefined \def \showCODEN     #1{\unskip}     \fi
\ifx \showDOI      \undefined \def \showDOI       #1{#1}\fi
\ifx \showISBNx    \undefined \def \showISBNx     #1{\unskip}     \fi
\ifx \showISBNxiii \undefined \def \showISBNxiii  #1{\unskip}     \fi
\ifx \showISSN     \undefined \def \showISSN      #1{\unskip}     \fi
\ifx \showLCCN     \undefined \def \showLCCN      #1{\unskip}     \fi
\ifx \shownote     \undefined \def \shownote      #1{#1}          \fi
\ifx \showarticletitle \undefined \def \showarticletitle #1{#1}   \fi
\ifx \showURL      \undefined \def \showURL       {\relax}        \fi
\providecommand\bibfield[2]{#2}
\providecommand\bibinfo[2]{#2}
\providecommand\natexlab[1]{#1}
\providecommand\showeprint[2][]{arXiv:#2}

\bibitem[Amit et~al\mbox{.}(2021)]%
        {amit2021segdiff}
\bibfield{author}{\bibinfo{person}{Tomer Amit}, \bibinfo{person}{Eliya
  Nachmani}, \bibinfo{person}{Tal Shaharabany}, {and} \bibinfo{person}{Lior
  Wolf}.} \bibinfo{year}{2021}\natexlab{}.
\newblock \showarticletitle{SegDiff: Image Segmentation with Diffusion
  Probabilistic Models}.
\newblock \bibinfo{journal}{\emph{CoRR}}  \bibinfo{volume}{abs/2112.00390}
  (\bibinfo{year}{2021}).
\newblock
\showeprint[arXiv]{2112.00390}
\urldef\tempurl%
\url{https://arxiv.org/abs/2112.00390}
\showURL{%
\tempurl}


\bibitem[Bayat(2023)]%
        {bayat2023a}
\bibfield{author}{\bibinfo{person}{Reza Bayat}.}
  \bibinfo{year}{2023}\natexlab{}.
\newblock \bibinfo{title}{A Study on Sample Diversity in Generative Models:
  {GAN}s vs. Diffusion Models}.
\newblock
\newblock
\urldef\tempurl%
\url{https://openreview.net/forum?id=BQpCuJoMykZ}
\showURL{%
\tempurl}


\bibitem[Bestuzheva et~al\mbox{.}(2023)]%
        {bestuzheva2021scip}
\bibfield{author}{\bibinfo{person}{Ksenia Bestuzheva}, \bibinfo{person}{Mathieu
  Besan{\c{c}}on}, \bibinfo{person}{Wei-Kun Chen}, \bibinfo{person}{Antonia
  Chmiela}, \bibinfo{person}{Tim Donkiewicz}, \bibinfo{person}{Jasper van
  Doornmalen}, \bibinfo{person}{Leon Eifler}, \bibinfo{person}{Oliver Gaul},
  \bibinfo{person}{Gerald Gamrath}, \bibinfo{person}{Ambros Gleixner},
  {et~al\mbox{.}}} \bibinfo{year}{2023}\natexlab{}.
\newblock \showarticletitle{The SCIP optimization suite 8.0}.
\newblock \bibinfo{journal}{\emph{{ACM} Trans. Math. Softw.}}
  \bibinfo{volume}{49}, \bibinfo{number}{2} (\bibinfo{year}{2023}),
  \bibinfo{pages}{22:1--22:21}.
\newblock
\urldef\tempurl%
\url{https://doi.org/10.1145/3585516}
\showDOI{\tempurl}


\bibitem[Dhariwal and Nichol(2021)]%
        {dhariwal2021diffusion}
\bibfield{author}{\bibinfo{person}{Prafulla Dhariwal} {and}
  \bibinfo{person}{Alexander Nichol}.} \bibinfo{year}{2021}\natexlab{}.
\newblock \showarticletitle{Diffusion models beat gans on image synthesis}.
\newblock \bibinfo{journal}{\emph{Advances in neural information processing
  systems}}  \bibinfo{volume}{34} (\bibinfo{year}{2021}),
  \bibinfo{pages}{8780--8794}.
\newblock


\bibitem[Gasse et~al\mbox{.}(2019)]%
        {gasse2019exact}
\bibfield{author}{\bibinfo{person}{Maxime Gasse}, \bibinfo{person}{Didier
  Ch{\'e}telat}, \bibinfo{person}{Nicola Ferroni}, \bibinfo{person}{Laurent
  Charlin}, {and} \bibinfo{person}{Andrea Lodi}.}
  \bibinfo{year}{2019}\natexlab{}.
\newblock \showarticletitle{Exact combinatorial optimization with graph
  convolutional neural networks}.
\newblock \bibinfo{journal}{\emph{Advances in Neural Information Processing
  Systems}}  \bibinfo{volume}{32} (\bibinfo{year}{2019}).
\newblock


\bibitem[Geifman and El-Yaniv(2019)]%
        {geifman2019selectivenet}
\bibfield{author}{\bibinfo{person}{Yonatan Geifman} {and} \bibinfo{person}{Ran
  El-Yaniv}.} \bibinfo{year}{2019}\natexlab{}.
\newblock \showarticletitle{Selectivenet: A deep neural network with an
  integrated reject option}. In \bibinfo{booktitle}{\emph{International
  conference on machine learning}}. PMLR, \bibinfo{pages}{2151--2159}.
\newblock


\bibitem[Gurobi~Optimization(2021)]%
        {gurobi2021gurobi}
\bibfield{author}{\bibinfo{person}{LLC Gurobi~Optimization}.}
  \bibinfo{year}{2021}\natexlab{}.
\newblock \bibinfo{title}{Gurobi optimizer reference manual}.
\newblock
\newblock


\bibitem[Han et~al\mbox{.}(2023)]%
        {han2023a}
\bibfield{author}{\bibinfo{person}{Qingyu Han}, \bibinfo{person}{Linxin Yang},
  \bibinfo{person}{Qian Chen}, \bibinfo{person}{Xiang Zhou},
  \bibinfo{person}{Dong Zhang}, \bibinfo{person}{Akang Wang},
  \bibinfo{person}{Ruoyu Sun}, {and} \bibinfo{person}{Xiaodong Luo}.}
  \bibinfo{year}{2023}\natexlab{}.
\newblock \showarticletitle{A {GNN}-Guided Predict-and-Search Framework for
  Mixed-Integer Linear Programming}. In \bibinfo{booktitle}{\emph{International
  Conference on Learning Representations}}.
\newblock
\urldef\tempurl%
\url{https://openreview.net/forum?id=pHMpgT5xWaE}
\showURL{%
\tempurl}


\bibitem[Ho et~al\mbox{.}(2020)]%
        {ho2020denoising}
\bibfield{author}{\bibinfo{person}{Jonathan Ho}, \bibinfo{person}{Ajay Jain},
  {and} \bibinfo{person}{Pieter Abbeel}.} \bibinfo{year}{2020}\natexlab{}.
\newblock \showarticletitle{Denoising diffusion probabilistic models}.
\newblock \bibinfo{journal}{\emph{Advances in neural information processing
  systems}}  \bibinfo{volume}{33} (\bibinfo{year}{2020}),
  \bibinfo{pages}{6840--6851}.
\newblock


\bibitem[Ho and Salimans(2022)]%
        {ho2022classifier}
\bibfield{author}{\bibinfo{person}{Jonathan Ho} {and} \bibinfo{person}{Tim
  Salimans}.} \bibinfo{year}{2022}\natexlab{}.
\newblock \showarticletitle{Classifier-Free Diffusion Guidance}.
\newblock \bibinfo{journal}{\emph{CoRR}}  \bibinfo{volume}{abs/2207.12598}
  (\bibinfo{year}{2022}).
\newblock
\urldef\tempurl%
\url{https://doi.org/10.48550/arXiv.2207.12598}
\showDOI{\tempurl}
\showeprint[arXiv]{2207.12598}


\bibitem[HUANG et~al\mbox{.}(2023)]%
        {huang2023accelerating}
\bibfield{author}{\bibinfo{person}{JUNWEI HUANG}, \bibinfo{person}{Zhiqing
  Sun}, {and} \bibinfo{person}{Yiming Yang}.} \bibinfo{year}{2023}\natexlab{}.
\newblock \showarticletitle{Accelerating Diffusion-based Combinatorial
  Optimization Solvers by Progressive Distillation}. In
  \bibinfo{booktitle}{\emph{ICML 2023 Workshop: Sampling and Optimization in
  Discrete Space}}.
\newblock
\urldef\tempurl%
\url{https://openreview.net/forum?id=AbMj31okE4}
\showURL{%
\tempurl}


\bibitem[Katoh and Ibaraki(1998)]%
        {katoh1998resource}
\bibfield{author}{\bibinfo{person}{Naoki Katoh} {and}
  \bibinfo{person}{Toshihide Ibaraki}.} \bibinfo{year}{1998}\natexlab{}.
\newblock \showarticletitle{Resource allocation problems}.
\newblock \bibinfo{journal}{\emph{Handbook of Combinatorial Optimization:
  Volume1--3}} (\bibinfo{year}{1998}), \bibinfo{pages}{905--1006}.
\newblock


\bibitem[Kelley(1960)]%
        {kelley1960cutting}
\bibfield{author}{\bibinfo{person}{James~E Kelley, Jr}.}
  \bibinfo{year}{1960}\natexlab{}.
\newblock \showarticletitle{The cutting-plane method for solving convex
  programs}.
\newblock \bibinfo{journal}{\emph{Journal of the society for Industrial and
  Applied Mathematics}} \bibinfo{volume}{8}, \bibinfo{number}{4}
  (\bibinfo{year}{1960}), \bibinfo{pages}{703--712}.
\newblock


\bibitem[Khalil et~al\mbox{.}(2022)]%
        {khalil2022mip}
\bibfield{author}{\bibinfo{person}{Elias~B Khalil},
  \bibinfo{person}{Christopher Morris}, {and} \bibinfo{person}{Andrea Lodi}.}
  \bibinfo{year}{2022}\natexlab{}.
\newblock \showarticletitle{Mip-gnn: A data-driven framework for guiding
  combinatorial solvers}. In \bibinfo{booktitle}{\emph{Proceedings of the AAAI
  Conference on Artificial Intelligence}}, Vol.~\bibinfo{volume}{36}.
  \bibinfo{pages}{10219--10227}.
\newblock


\bibitem[Kingma and Ba(2015)]%
        {kingma2014adam}
\bibfield{author}{\bibinfo{person}{Diederik~P. Kingma} {and}
  \bibinfo{person}{Jimmy Ba}.} \bibinfo{year}{2015}\natexlab{}.
\newblock \showarticletitle{Adam: {A} Method for Stochastic Optimization}. In
  \bibinfo{booktitle}{\emph{3rd International Conference on Learning
  Representations, {ICLR} 2015, San Diego, CA, USA, May 7-9, 2015, Conference
  Track Proceedings}}, \bibfield{editor}{\bibinfo{person}{Yoshua Bengio} {and}
  \bibinfo{person}{Yann LeCun}} (Eds.).
\newblock
\urldef\tempurl%
\url{http://arxiv.org/abs/1412.6980}
\showURL{%
\tempurl}


\bibitem[Lawler and Wood(1966)]%
        {lawler1966branch}
\bibfield{author}{\bibinfo{person}{Eugene~L Lawler} {and}
  \bibinfo{person}{David~E Wood}.} \bibinfo{year}{1966}\natexlab{}.
\newblock \showarticletitle{Branch-and-bound methods: A survey}.
\newblock \bibinfo{journal}{\emph{Operations research}} \bibinfo{volume}{14},
  \bibinfo{number}{4} (\bibinfo{year}{1966}), \bibinfo{pages}{699--719}.
\newblock


\bibitem[Li et~al\mbox{.}(2023)]%
        {li2023qdiffusion}
\bibfield{author}{\bibinfo{person}{Xiuyu Li}, \bibinfo{person}{Yijiang Liu},
  \bibinfo{person}{Long Lian}, \bibinfo{person}{Huanrui Yang},
  \bibinfo{person}{Zhen Dong}, \bibinfo{person}{Daniel Kang},
  \bibinfo{person}{Shanghang Zhang}, {and} \bibinfo{person}{Kurt Keutzer}.}
  \bibinfo{year}{2023}\natexlab{}.
\newblock \showarticletitle{Q-Diffusion: Quantizing Diffusion Models}. In
  \bibinfo{booktitle}{\emph{Proceedings of the IEEE/CVF International
  Conference on Computer Vision (ICCV)}}. \bibinfo{pages}{17535--17545}.
\newblock


\bibitem[Loshchilov and Hutter(2019)]%
        {LoshchilovH19}
\bibfield{author}{\bibinfo{person}{Ilya Loshchilov} {and}
  \bibinfo{person}{Frank Hutter}.} \bibinfo{year}{2019}\natexlab{}.
\newblock \showarticletitle{Decoupled Weight Decay Regularization}. In
  \bibinfo{booktitle}{\emph{7th International Conference on Learning
  Representations, {ICLR} 2019, New Orleans, LA, USA, May 6-9, 2019}}.
  \bibinfo{publisher}{OpenReview.net}.
\newblock
\urldef\tempurl%
\url{https://openreview.net/forum?id=Bkg6RiCqY7}
\showURL{%
\tempurl}


\bibitem[Maher et~al\mbox{.}(2017)]%
        {maher2017scip}
\bibfield{author}{\bibinfo{person}{Stephen~J Maher}, \bibinfo{person}{Tobias
  Fischer}, \bibinfo{person}{Tristan Gally}, \bibinfo{person}{Gerald Gamrath},
  \bibinfo{person}{Ambros Gleixner}, \bibinfo{person}{Robert~Lion Gottwald},
  \bibinfo{person}{Gregor Hendel}, \bibinfo{person}{Thorsten Koch},
  \bibinfo{person}{Marco L{\"u}bbecke}, \bibinfo{person}{Matthias
  Miltenberger}, {et~al\mbox{.}}} \bibinfo{year}{2017}\natexlab{}.
\newblock \showarticletitle{The SCIP optimization suite 4.0}.
\newblock  (\bibinfo{year}{2017}).
\newblock


\bibitem[Maz{\'e} and Ahmed(2023)]%
        {maze2023diffusion}
\bibfield{author}{\bibinfo{person}{Fran{\c{c}}ois Maz{\'e}} {and}
  \bibinfo{person}{Faez Ahmed}.} \bibinfo{year}{2023}\natexlab{}.
\newblock \showarticletitle{Diffusion models beat gans on topology
  optimization}. In \bibinfo{booktitle}{\emph{Proceedings of the AAAI
  Conference on Artificial Intelligence (AAAI), Washington, DC}}.
\newblock


\bibitem[Nair et~al\mbox{.}(2020)]%
        {nair2020solving}
\bibfield{author}{\bibinfo{person}{Vinod Nair}, \bibinfo{person}{Sergey
  Bartunov}, \bibinfo{person}{Felix Gimeno}, \bibinfo{person}{Ingrid
  Von~Glehn}, \bibinfo{person}{Pawel Lichocki}, \bibinfo{person}{Ivan Lobov},
  \bibinfo{person}{Brendan O'Donoghue}, \bibinfo{person}{Nicolas Sonnerat},
  \bibinfo{person}{Christian Tjandraatmadja}, \bibinfo{person}{Pengming Wang},
  {et~al\mbox{.}}} \bibinfo{year}{2020}\natexlab{}.
\newblock \showarticletitle{Solving mixed integer programs using neural
  networks}.
\newblock \bibinfo{journal}{\emph{arXiv preprint arXiv:2012.13349}}
  (\bibinfo{year}{2020}).
\newblock


\bibitem[Nichol and Dhariwal(2021)]%
        {nichol2021improved}
\bibfield{author}{\bibinfo{person}{Alexander~Quinn Nichol} {and}
  \bibinfo{person}{Prafulla Dhariwal}.} \bibinfo{year}{2021}\natexlab{}.
\newblock \showarticletitle{Improved denoising diffusion probabilistic models}.
  In \bibinfo{booktitle}{\emph{International Conference on Machine Learning}}.
  PMLR, \bibinfo{pages}{8162--8171}.
\newblock


\bibitem[Pantelides et~al\mbox{.}(1995)]%
        {pantelides1995short}
\bibfield{author}{\bibinfo{person}{CC Pantelides}, \bibinfo{person}{MJ Realff},
  {and} \bibinfo{person}{N Shah}.} \bibinfo{year}{1995}\natexlab{}.
\newblock \showarticletitle{Short-term scheduling of pipeless batch plants}.
\newblock \bibinfo{journal}{\emph{Chemical engineering research \& design}}
  \bibinfo{volume}{73}, \bibinfo{number}{4} (\bibinfo{year}{1995}),
  \bibinfo{pages}{431--444}.
\newblock


\bibitem[Pisinger and Ropke(2019)]%
        {pisinger2019large}
\bibfield{author}{\bibinfo{person}{David Pisinger} {and}
  \bibinfo{person}{Stefan Ropke}.} \bibinfo{year}{2019}\natexlab{}.
\newblock \showarticletitle{Large neighborhood search}.
\newblock \bibinfo{journal}{\emph{Handbook of metaheuristics}}
  (\bibinfo{year}{2019}), \bibinfo{pages}{99--127}.
\newblock


\bibitem[Pochet and Wolsey(2006)]%
        {pochet2006}
\bibfield{author}{\bibinfo{person}{Yves Pochet} {and}
  \bibinfo{person}{Laurence~A Wolsey}.} \bibinfo{year}{2006}\natexlab{}.
\newblock \bibinfo{booktitle}{\emph{Production planning by mixed integer
  programming}}. Vol.~\bibinfo{volume}{149}.
\newblock \bibinfo{publisher}{Springer}.
\newblock


\bibitem[Prouvost et~al\mbox{.}(2020)]%
        {prouvost2020ecole}
\bibfield{author}{\bibinfo{person}{Antoine Prouvost}, \bibinfo{person}{Justin
  Dumouchelle}, \bibinfo{person}{Lara Scavuzzo}, \bibinfo{person}{Maxime
  Gasse}, \bibinfo{person}{Didier Ch{\'e}telat}, {and} \bibinfo{person}{Andrea
  Lodi}.} \bibinfo{year}{2020}\natexlab{}.
\newblock \showarticletitle{Ecole: A gym-like library for machine learning in
  combinatorial optimization solvers}.
\newblock \bibinfo{journal}{\emph{arXiv preprint arXiv:2011.06069}}
  (\bibinfo{year}{2020}).
\newblock


\bibitem[Radford et~al\mbox{.}(2021)]%
        {radford2021learning}
\bibfield{author}{\bibinfo{person}{Alec Radford}, \bibinfo{person}{Jong~Wook
  Kim}, \bibinfo{person}{Chris Hallacy}, \bibinfo{person}{Aditya Ramesh},
  \bibinfo{person}{Gabriel Goh}, \bibinfo{person}{Sandhini Agarwal},
  \bibinfo{person}{Girish Sastry}, \bibinfo{person}{Amanda Askell},
  \bibinfo{person}{Pamela Mishkin}, \bibinfo{person}{Jack Clark},
  {et~al\mbox{.}}} \bibinfo{year}{2021}\natexlab{}.
\newblock \showarticletitle{Learning transferable visual models from natural
  language supervision}. In \bibinfo{booktitle}{\emph{International conference
  on machine learning}}. PMLR, \bibinfo{pages}{8748--8763}.
\newblock


\bibitem[Ramesh et~al\mbox{.}(2022)]%
        {ramesh2022hierarchical}
\bibfield{author}{\bibinfo{person}{Aditya Ramesh}, \bibinfo{person}{Prafulla
  Dhariwal}, \bibinfo{person}{Alex Nichol}, \bibinfo{person}{Casey Chu}, {and}
  \bibinfo{person}{Mark Chen}.} \bibinfo{year}{2022}\natexlab{}.
\newblock \showarticletitle{Hierarchical text-conditional image generation with
  clip latents}.
\newblock \bibinfo{journal}{\emph{arXiv preprint arXiv:2204.06125}}
  \bibinfo{volume}{1}, \bibinfo{number}{2} (\bibinfo{year}{2022}),
  \bibinfo{pages}{3}.
\newblock


\bibitem[Sawik(2011)]%
        {sawik2011}
\bibfield{author}{\bibinfo{person}{Tadeusz Sawik}.}
  \bibinfo{year}{2011}\natexlab{}.
\newblock \bibinfo{booktitle}{\emph{Scheduling in supply chains using mixed
  integer programming}}.
\newblock \bibinfo{publisher}{John Wiley \& Sons}.
\newblock


\bibitem[Silver et~al\mbox{.}(1998)]%
        {silver1998inventory}
\bibfield{author}{\bibinfo{person}{Edward~Allen Silver},
  \bibinfo{person}{David~F Pyke}, \bibinfo{person}{Rein Peterson},
  {et~al\mbox{.}}} \bibinfo{year}{1998}\natexlab{}.
\newblock \bibinfo{booktitle}{\emph{Inventory management and production
  planning and scheduling}}. Vol.~\bibinfo{volume}{3}.
\newblock \bibinfo{publisher}{Wiley New York}.
\newblock


\bibitem[Sohl-Dickstein et~al\mbox{.}(2015)]%
        {sohl2015deep}
\bibfield{author}{\bibinfo{person}{Jascha Sohl-Dickstein},
  \bibinfo{person}{Eric Weiss}, \bibinfo{person}{Niru Maheswaranathan}, {and}
  \bibinfo{person}{Surya Ganguli}.} \bibinfo{year}{2015}\natexlab{}.
\newblock \showarticletitle{Deep unsupervised learning using nonequilibrium
  thermodynamics}. In \bibinfo{booktitle}{\emph{International conference on
  machine learning}}. PMLR, \bibinfo{pages}{2256--2265}.
\newblock


\bibitem[Sohn et~al\mbox{.}(2015)]%
        {sohn2015learning}
\bibfield{author}{\bibinfo{person}{Kihyuk Sohn}, \bibinfo{person}{Honglak Lee},
  {and} \bibinfo{person}{Xinchen Yan}.} \bibinfo{year}{2015}\natexlab{}.
\newblock \showarticletitle{Learning structured output representation using
  deep conditional generative models}.
\newblock \bibinfo{journal}{\emph{Advances in neural information processing
  systems}}  \bibinfo{volume}{28} (\bibinfo{year}{2015}).
\newblock


\bibitem[Song et~al\mbox{.}(2021a)]%
        {song2020denoising}
\bibfield{author}{\bibinfo{person}{Jiaming Song}, \bibinfo{person}{Chenlin
  Meng}, {and} \bibinfo{person}{Stefano Ermon}.}
  \bibinfo{year}{2021}\natexlab{a}.
\newblock \showarticletitle{Denoising Diffusion Implicit Models}. In
  \bibinfo{booktitle}{\emph{9th International Conference on Learning
  Representations, {ICLR} 2021, Virtual Event, Austria, May 3-7, 2021}}.
  \bibinfo{publisher}{OpenReview.net}.
\newblock
\urldef\tempurl%
\url{https://openreview.net/forum?id=St1giarCHLP}
\showURL{%
\tempurl}


\bibitem[Song et~al\mbox{.}(2021b)]%
        {song2020score}
\bibfield{author}{\bibinfo{person}{Yang Song}, \bibinfo{person}{Jascha
  Sohl{-}Dickstein}, \bibinfo{person}{Diederik~P. Kingma},
  \bibinfo{person}{Abhishek Kumar}, \bibinfo{person}{Stefano Ermon}, {and}
  \bibinfo{person}{Ben Poole}.} \bibinfo{year}{2021}\natexlab{b}.
\newblock \showarticletitle{Score-Based Generative Modeling through Stochastic
  Differential Equations}. In \bibinfo{booktitle}{\emph{9th International
  Conference on Learning Representations, {ICLR} 2021, Virtual Event, Austria,
  May 3-7, 2021}}. \bibinfo{publisher}{OpenReview.net}.
\newblock
\urldef\tempurl%
\url{https://openreview.net/forum?id=PxTIG12RRHS}
\showURL{%
\tempurl}


\bibitem[Toth and Vigo(2002)]%
        {toth2002vehicle}
\bibfield{author}{\bibinfo{person}{Paolo Toth} {and} \bibinfo{person}{Daniele
  Vigo}.} \bibinfo{year}{2002}\natexlab{}.
\newblock \bibinfo{booktitle}{\emph{The vehicle routing problem}}.
\newblock \bibinfo{publisher}{SIAM}.
\newblock


\bibitem[Wolsey(1998)]%
        {wolsey1998}
\bibfield{author}{\bibinfo{person}{L.A. Wolsey}.}
  \bibinfo{year}{1998}\natexlab{}.
\newblock \bibinfo{booktitle}{\emph{Integer Programming}}.
\newblock \bibinfo{publisher}{Wiley}.
\newblock
\showISBNx{9780471283669}
\showLCCN{98007296}
\urldef\tempurl%
\url{https://books.google.co.uk/books?id=x7RvQgAACAAJ}
\showURL{%
\tempurl}


\bibitem[Yang et~al\mbox{.}(2022)]%
        {yang2022diffusion}
\bibfield{author}{\bibinfo{person}{Ling Yang}, \bibinfo{person}{Zhilong Zhang},
  \bibinfo{person}{Yang Song}, \bibinfo{person}{Shenda Hong},
  \bibinfo{person}{Runsheng Xu}, \bibinfo{person}{Yue Zhao},
  \bibinfo{person}{Yingxia Shao}, \bibinfo{person}{Wentao Zhang},
  \bibinfo{person}{Ming{-}Hsuan Yang}, {and} \bibinfo{person}{Bin Cui}.}
  \bibinfo{year}{2022}\natexlab{}.
\newblock \showarticletitle{Diffusion Models: {A} Comprehensive Survey of
  Methods and Applications}.
\newblock \bibinfo{journal}{\emph{CoRR}}  \bibinfo{volume}{abs/2209.00796}
  (\bibinfo{year}{2022}).
\newblock
\urldef\tempurl%
\url{https://doi.org/10.48550/arXiv.2209.00796}
\showDOI{\tempurl}
\showeprint[arXiv]{2209.00796}


\bibitem[Yoon(2022)]%
        {yoon2022confidence}
\bibfield{author}{\bibinfo{person}{Taehyun Yoon}.}
  \bibinfo{year}{2022}\natexlab{}.
\newblock \showarticletitle{Confidence Threshold Neural Diving}.
\newblock \bibinfo{journal}{\emph{CoRR}}  \bibinfo{volume}{abs/2202.07506}
  (\bibinfo{year}{2022}).
\newblock
\showeprint[arXiv]{2202.07506}
\urldef\tempurl%
\url{https://arxiv.org/abs/2202.07506}
\showURL{%
\tempurl}


\end{thebibliography}

\appendix

\section{Appendix}

\subsection{CISP Algorithm}
\label{app:cisp}
The study conducted by \cite{radford2021learning} underscores the substantial efficacy of contrasting pre-training in capturing multi-modal data, with particular emphasis on its application in the text-to-image transfer domain. Drawing inspiration from this seminal work, we introduce the CISP algorithm. The primary objective of CISP is to facilitate the learning of the IP Encoder $\mathbf{E}_{I}$ and the Solution Encoder $\mathbf{E}_{\mathbf{X}}$, as illustrated in Algorithm \ref{alg:cisp}. Within the scope of our investigation, we work with a mini-batch of data comprising instances denoted as $I$ and their corresponding solutions denoted as $\mathbf{X}$. The batch's bipartite graph representing instances $I$ is denoted as $\mathbf{G}$. We use $\mathbf{z}_{I}$ and $\mathbf{z}_{\mathbf{X}}$ to denote the embeddings of instances and solutions, respectively. Notably, both $\mathbf{z}_{I}$ and $\mathbf{z}_{\mathbf{X}}$ possess identical dimensions, enabling us to compute their cosine similarity. Furthermore, within the mini-batch,$\mathbf{z}_{I,j}$ and $\mathbf{z}_{\mathbf{X},k}$ denote to the $j$th sample in $\mathbf{z}_{I}$ and the $k$th sample in $\mathbf{z}_{\mathbf{X}}$, respectively. Within this conceptual framework, we leverage the matrix $\mathbf{s}\in \mathbb{R}^{N\times N}$ to represent the similarity between $N$ instances and $N$ solutions. Each element $\mathbf{s}_{j,k}$, where $j,k\in \{1,...,N\}$, corresponds to the logit employed in computation of the symmetric cross-entropy loss.
\begin{algorithm}[ht]
\caption{Contrastive IP-Solution Pre-Training (CISP)}
\label{alg:cisp}
\begin{algorithmic}[1] 
\ENSURE The mini-batch size $N$, the mini-batch bipartite graph representations of IP instance set $I$, denoted by $\mathbf{G}$,  and corresponding mini-batch solutions $\mathbf{X}$ \\
\REQUIRE IP Encoder $\mathbf{E}_{I}$, Solution Encoder $\mathbf{E}_{\mathbf{X}}$, temperature parameter $\tau$
    \STATE Get IP and solution embeddings $\mathbf{z}_{I}, \mathbf{z}_{\mathbf{X}} = \mathbf{E}_{I}(\mathbf{G}), \mathbf{E}_{\mathbf{X}}(\mathbf{X})$ {\color{gray}\small \textit{   // $N\times n\times d$, where $n$ is the padding length of variables and $d$ is the embedding size.}}
    \FOR{$j \in \{1,2,..., N\}$ and $k \in \{1,2,..., N\}$} 
        \STATE Flatten $\mathbf{z}_{I,j}$ and $\mathbf{z}_{\mathbf{X},k}$ into vectors $\bar{\mathbf{z}}_{I,j}$ and $\bar{\mathbf{z}}_{\mathbf{X},k}$
        \STATE $\mathbf{s}_{j,k} = e^{\tau} \cdot \bar{\mathbf{z}}_{I,j}^T\bar{\mathbf{z}}_{\mathbf{X},k}/(\|\bar{\mathbf{z}}_{I,j}\| \|\bar{\mathbf{z}}_{\mathbf{X},k}\|)$ {\color{gray}\small \textit{   // compute similarity for IP and solution embeddings }}
    \ENDFOR
    \STATE Set labels $\mathbf{y} = (1,2,..., N)$ 
    \STATE Compute cross-entropy loss $\mathcal{L}_I$ by utilizing $\mathbf{s}_{j,*}$ and $\mathbf{y}$.
    \STATE  Compute cross-entropy loss $\mathcal{L}_{\mathbf{X}}$ by utilizing $\mathbf{s}_{*,k}$ and $\mathbf{y}$. 
    \STATE Compute the symmetric loss $\mathcal{L} =(\mathcal{L}_I+ \mathcal{L}_{\mathbf{X}})/2$
    \RETURN $\mathcal{L}$
\end{algorithmic}
\end{algorithm}

\subsection{Feature Descriptions For Variables Nodes, Constraint Nodes And Edges}
\label{app:features}
In Table \ref{tab:feature_description}, we provide a description of the features that are extracted using the \textit{Ecole} library \citep{prouvost2020ecole} and used as IP bipartite graph representations for training the GCN model.
\begin{table*}[ht]
\caption{Description of the variable, constraint and edge features in our bipartite graph representations.}
    \begin{tabular}{lll}
    \toprule
                          & \multicolumn{1}{c}{\small Feature} & \multicolumn{1}{c}{\small Description}                                         \\ \midrule
    \multirow{11}{*}{\small Variable}  & \small type                        & \small Type(binary, integer, impl. integer, continuous) as a one-hot encoding. \\ \cline{2-3} 
         & \small coef         & \small Objective coefficient, normalized.                   \\ \cline{2-3} 
         & \small has\_lb       & \small Lower bound indicator.                               \\ \cline{2-3} 
         & \small has\_ub       & \small Upper bound indicator.                               \\ \cline{2-3} 
         & \small sol\_is\_at\_lb & \small Solution value equals lower bound.                   \\ \cline{2-3} 
         & \small sol\_is\_at\_ub & \small Solution value equals upper bound.                   \\ \cline{2-3} 
         & \small sol\_frac     & \small Solution value fractionality.                        \\ \cline{2-3} 
         & \small basis\_status   & \small Simplex basis status(lower, basic, upper, zero) as a one-hot encoding.  
                \\ \cline{2-3} 
         & \small reduced\_cost & \small Reduced cost, normalized.                            \\ \cline{2-3} 
         & \small age          & \small LP age, normalized.                                  \\ \cline{2-3} 
         & \small sol\_val      & \small Solution value.                                      \\ \midrule
        \multirow{5}{*}{\small Constraint} & \small obj\_cos\_sim                 & \small Cosine similarity with objective.                                       \\ \cline{2-3} 
         & \small bias         & \small Bias value, normalized with constraint coefficients. \\ \cline{2-3} 
         & \small is\_tight     & \small Tightness indicator in LP solution.                  \\ \cline{2-3} 
         & \small dualsol\_val  & \small Dual solution value, normalized.                     \\ \cline{2-3} 
         & \small age          & \small LP age, normalized with total number of LPs.         \\ \midrule
    \small Edge & \small coef         & \small Constraint coefficient, normalized per constraint.   \\ \bottomrule
    \end{tabular}
    \label{tab:feature_description}
\end{table*}

\subsection{Datasets and Baselines}
\label{app:data_baselines}

\subsubsection{Datasets.} 
For all four datasets, we randomly generate 1000 instances (800 for training, 100 for validation and 100 for testing). 
Table~\ref{tab:datasize} summarizes the numbers of constraints, variables and problem type of each dataset.
We then collect feasible solutions and their objective values for each instance by running the Gurobi~\citep{gurobi2021gurobi} or SCIP~\cite{bestuzheva2021scip}, 
where the time limit is set to 1000s for each instance. 
For those instances with a large number of feasible solutions, 
we only keep 500 best solutions. 
We adopt the same features exacted via the Ecole library as in~\citep{gasse2019exact} and exclude those related to feasible solutions. The specific features are shown in Appendix \ref{app:features}.

\subsubsection{Baselines.} 
We compared our method with the following baselines: 
\begin{itemize}
\item \textit{SCIP}~\citep{bestuzheva2021scip} (an open source solver): SCIP is currently one of the fastest non-commercial solvers for mixed integer programming (MIP) and mixed integer nonlinear programming (MINLP). Here, we focus on comparing the quality of initial solutions with SCIP. Instead of relying on the first feasible solution generated by SCIP, which are often of low quality due to the use of trivial heuristics, we employ the first solution produced via non-trivial heuristic algorithms during the solving process of SCIP \citep{bestuzheva2021scip} (i.e. the first feasible solution after the pre-solving stage of SCIP).
\item \textit{Gurobi}~\citep{gurobi2021gurobi} (the powerful commercial solver): Gurobi is a highly efficient commercial mathematical optimization solver. As our focus is on comparing feasible solutions, we consider the best solutions obtained through Gurobi's default heuristic algorithms as a benchmark.
\item \textit{Neural Diving (ND)}~\citep{nair2020solving}: Neural Diving adopts a solution prediction approach, training a GCN to predict the value of each variable. It then incorporates SelectiveNet \citep{geifman2019selectivenet} to generate partial solutions with a predefined coverage threshold $C$ (e.g., a coverage threshold $C=0.2$ means that the expectation of the number of assigned variables by the neural network is 20\%). This threshold is typically set to be low ($<$0.5) to ensure the feasibility of partial solutions.
In our experiments, to evaluate the feasibility of solutions, we set two different coverage levels, namely low coverage (which usually indicates higher feasibility) and high coverage (which usually indicates lower feasibility), for each dataset. We then compare the feasibility of the partial solutions with the complete solutions generated by our methods. To assess the quality of solutions, we employ the Completesol heuristic in SCIP (Algorithm 4 in \cite{maher2017scip}) to enhance the partial solutions. This heuristic involves solving auxiliary IP instances by fixing the variables from the partial solutions. By utilizing this heuristic, we can obtain more complete solutions and evaluate their quality.
\item \textit{Predict-and-search algorithm (PS)}~\citep{han2023a}: The PS algorithm, similar to Neural Diving, utilizes graph neural networks to predict the value of each variable. It then searches for the best feasible solution within a trust region constructed by the predicted partial solutions. This method requires setting parameters $(k_0, k_1)$ to represent the numbers of 0's and 1's in a partial solution, and $\Delta$ to define the range of the neighborhood region of the partial solution. To search a high-quality feasible solution, this method adds neighborhood constraints to origin instance, which produces modified IP instance. Therefore, an IP solver such as SCIP or Gurobi is required to solve the modified instance and obtain feasible solutions. In our experiments, we use Gurobi as the solver and control the parameters $\Delta$ to ensure that the modified instance is 100\% feasible. We considered the best heuristic solutions from the modified instance found by Gurobi as our baseline.
\end{itemize}

\begin{table}[ht]
    \centering
     \caption{Instance size of each dataset}
    \begin{tabular}{cccc}
        \toprule
        \small Dataset & \small Constraints &\small Variables & \small Problem Type\\
        \midrule
        \small SC & \small 1000 & \small 2000 & \small minimize\\
        \small CF & \small 5051 & \small 5050 & \small minimize\\
        \small CA & \small 786 & \small 1500 & \small maximize\\
        \small IS & \small 6396 & \small 1500 & \small maximize\\
        \bottomrule
    \end{tabular}
    \label{tab:datasize}
\end{table}

\subsection{Training Details}
\label{app:training}

We trained the CISP and diffusion model on four IP datasets. Each dataset contained 800 training instances, with 500 solutions collected for each instance. In each batch, we sampled $64$ instances, and for each instance, We sample one solution from 500 solutions in proportion to the probability of the objective value as a possible label. This implies that solutions with better objective values had a higher probability of being sampled. We iterated through all instances (with one solution per instance) in each epoch.

For the Solution Encoder, we utilized a single transformer encode layer with a width of 128. The IP encoder adopted the architecture described in \cite{nair2020solving}, using GCN to obtain embeddings for all variables as IP embeddings. Both models transformed the features of the solution and IP into latent variables with a dimension of 128, enabling convenient computation of cosine similarity in CISP. The CISP was trained using the AdamW Optimizer \citep{LoshchilovH19}. We employed a decreasing learning rate strategy, starting with a learning rate of 0.001 and linearly decaying it by a factor of 0.9 every 100 epochs until reaching 800 epochs. The model training was performed with a batch size of 64.

For the diffusion model, we utilized a single-layer Transformer encoder with a width of 128 to predict $\bzx$ and adjusted the number of time steps to 1000. The forward process variances were set as constants, increasing linearly from $\beta_1 = 10^{-4}$ to $\beta_T = 0.02$, following the default setting of DDPM \citep{ho2020denoising}. The solution decoder model was jointly trained with the diffusion model and consisted of two Transformer encode layers with a width of 128. The loss function was defined as the sum of the diffusion loss, decoder loss, and the penalty for violating constraints, as shown in \eqref{loss:total}. Here, $\lambda$ is set to be the number of variables in the instances from the training set, excluding the IS dataset, where $\lambda = 0$. We trained diffusion and decoder model for 100 epochs with batch size of $32$ via Adam Optimizer \citep{kingma2014adam}.

\subsection{Training and Inference Time}
\label{sec:train_inf_time}
In this section, we report the training time (including CISP pretraining and Diffusion model) for each dataset, 
which takes 100 epochs for both CISP and Diffusion model to converge. 
Additionally, we provide the total inference time for sampling 3000 solutions by using IP Guided DDIM and DDPM.
From Table \ref{tab:time}, we observe that our method requires a reasonable amount of time for model training. 
During the inference phase, IP Guided DDIM demonstrates faster performance compared to IP Guided DDPM with average time of 0.46s-1.68s for sampling each solution. 
Moreover, as shown in the experiment results from Section \ref{sec:experiments}, 
IP Guided DDIM also achieves better performance than DDPM, making it suitable for practical applications.
\begin{table}[ht]
    \centering
    \caption{Total training time and inference time for sampling 3000 solutions for each dataset}
    \resizebox{\linewidth}{!}{
    \begin{tabular}{cccccc}
            \toprule
            \small Dataset & \small Training (CISP + Diffusion) &\small IP Guided DDIM   & \small IP Guided DDPM  \\
            \midrule
            \small SC & \small 24.4m & \small 37.5m & \small  374m\\
            \small CF & \small 71.7m & \small 84m & \small 805m \\
        \small CA & \small 9.3m & \small  23m & \small 233.5m \\
        \small IS & \small 11.1m & \small 23m & \small 234m \\
            \bottomrule
        \end{tabular}}
    \label{tab:time}
\end{table}

\subsection{Hyperparameters for IP Guided Sampling}
\label{app:sampling-parameter}
During the sampling process, we configured the number of steps to be 1000 for IP Guided Diffusion sampling (IP Guided DDPM) and 100 for Non-Markovian IP Guided Diffusion sampling (IP Guided DDIM). In Table \ref{tab:hyperparameter_settings}, we provide the specific values for the gradient scale $s$ and leverage factor $\gamma$.
 
\begin{table}[ht]
\centering
\caption{$s$ and $\gamma$ settings in different dataset}
\begin{tabular}{lcccc}
\toprule
\multicolumn{1}{c}{\multirow{2}{*}{\small dataset}} & \multicolumn{2}{c}{\small IP Guided DDIM} & \multicolumn{2}{c}{\small IP Guided DDPM} \\ \cline{2-5} 
\multicolumn{1}{c}{}     & \small $s$  & \small $\gamma$   & \small $s$  & \small $\gamma$  \\ \midrule
\small SC (2000)  & \small 100,000 & \small 0.9 & \small 15,000  & \small 0.1 \\
\small SC (3000) & \small 150,000 & \small 0.9 & \small 22,500  & \small 0.1 \\
\small SC (4000)  & \small 200,000 & \small 0.9 & \small 30,000  & \small 0.1 \\
\small CF         & \small 1,000   & \small 0.7 & \small 500,000 & \small 0.1 \\
\small CA         & \small 20,000  & \small 0.7 & \small 10,000  & \small 0.3 \\
\small IS         & \small 20,000  & \small 0.5 & \small 10,000  & \small 0.1 \\ \bottomrule
\end{tabular}
\label{tab:hyperparameter_settings}
\end{table}

\end{document}